\documentclass[final,epsfig,amsfont]{article}
\usepackage{latexsym}
\usepackage[dvips]{graphicx}
\newtheorem{theorem}{Theorem}[section]
\newtheorem{lemma}[theorem]{Lemma}
\newtheorem{proposition}[theorem]{Proposition}

\newtheorem{remark}[theorem]{Remark}

\newtheorem{definition} {Definition} [section]

\def\Cal{\mathcal}

\def\ccb{{\Cal B}}

\def\ck{{\Cal K}}
\def\cull{{Conv}}

\def\enpr{\quad \vrule height .9ex width .8ex depth -.1ex}

\def\bar{\overline}

\def\bre{\mathbf R}
\def\bna{\mathbf N}

\def\llim{\mathop{\longrightarrow}}

\def\s,{\quad $\,$}

\def\smad{\smallskip Proof.\ }

\def\s,{\quad $\,$}

\def\witi{\widetilde}

\def\disp{\displaystyle}

\def\bfu  {{\bf U}}
\def\bfup  {{\bf U}_{\Psi}}
\def\bfuf  {{\bf U}_{\Phi}}
\def\wii{\Cal I}
\def\wiib{\overline{\Cal I}}

\begin{document}

\title{Attractors of Iterated Function Systems with uncountably many maps}
\author{
Giorgio Mantica \\
Center for Non-linear and Complex Systems,\\
Dipartimento di Scienze ed Alta Tecnologia, \\ Universit\`a dell' Insubria,
via Valleggio 11, 22100
Como, Italy. \\ Also at I.N.F.N. sezione di Milano and CNISM unit\`a di Como.\\
\\
Roberto Peirone \\
Dipartimento di Matematica, \\Universit\`a di Roma Tor Vergata, \\ via della Ricerca Scientifica, 00133 Roma, Italy.
}
\date{}
\maketitle

\begin{abstract}
We study the topological properties of attractors of Iterated Function Systems (I.F.S.) on the real line, consisting of affine maps of homogeneous contraction ratio. These maps define what we call a {\em second generation} I.F.S.: they are uncountably many and the set of their fixed points is a Cantor set. We prove that when this latter either is the attractor of a finite, non-singular, hyperbolic, I.F.S. (of {\em first generation}), or it possesses a particular {\em dissection} property, the attractor of the second generation I.F.S. consists of finitely many closed intervals.
\end{abstract}


{\em keywords}
Cantor sets, Iterated Function Systems, Second Generation I.F.S., Attractors.

\section{Introduction and discussion of the main results}

Let $\Psi$ be a set of contractive transformations on $\bre$. Also, let the operator $\bfup$ be defined by
\begin{equation}
\label{eq-u1}
 \bfup(A)=\overline{{\bigcup\limits_{\psi\in\Psi}  \psi(A)}}
\end{equation}
for every $A\subseteq\bre$, where the bar denotes topological closure.
Let $\ck$ be the set of nonempty compact subsets of $\bre$, endowed with the Hausdorff metric \cite{ba2}. In this setting, $\bfup$ is a contractive operator on $\ck$ and $K =K_\Psi$
is the unique element that solves the equation

\begin{equation}
K=\bfup\big(K\big).
\label{eqno(1)}
\end{equation}

$K_\Psi$ is termed the attractor of the Iterated Functions System $\Psi$ \cite{papmor,hut,dem,hata,duvall}.
Following standard terminology, an I.F.S. consisting of contractive maps is called {\em hyperbolic}.

When the cardinality of $\Psi$ is finite, $K_\Psi$ may be an interval, a countable union of intervals, and a Cantor set. Attractors of I.F.S. with
countable sets of maps have been considered in \cite{moran,urba,fernau,molter,hille}. Their study requires the introduction of the topological closure of the r.h.s. of eq. (\ref{eq-u1}), at difference with the case of finitely many transformations.

In this paper we will study a more general case composed of a continuous set of maps, but we will restrict ourselves in two ways. First, following Elton and Yan \cite{elton} we will consider {\em homogeneous} affine maps. Secondly, as in  \cite{nalgo2,arxiv,intj} these maps will be structured as a {\em second generation I.F.S.}.

Precisely, we start from a {\em first generation I.F.S.} $\Psi$ consisting of a finite number $M$ of maps,
\begin{equation}
\Psi:=\{\psi_i: i=1,\ldots,M\}.
\label{eq-map00}
\end{equation}
We assume that every map $\psi_i$ is $C^2$ and that there exist  constants $\delta$ and $\sigma$ such that for any $i=1,\ldots,M$
\begin{equation}
   0 < \sigma \leq |\psi'_i(x)| \leq  \delta < 1
\label{eq-contra}
\end{equation}
for all $x$ in a suitable interval $I$ including all fixed points of the maps $\psi_i$. We also suppose that $\psi_i(I) \subseteq I$ for any $i$, and that there exist at least two different fixed points.
All the above conditions define what we call a {\em non-singular, hyperbolic I.F.S.}
As before, $K_\Psi$ denotes the attractor of this I.F.S.

Next, consider a new set of affine maps, of equal contraction ratio $0 <\alpha<1$, and where $\beta \in \bre $:
\begin{equation}
\phi(\beta;x)=\alpha (x-\beta)+\beta=\alpha x +(1-\alpha) \beta.
\label{eq-map1}
\end{equation}
We will also use the notation $\phi_\beta(\cdot) = \phi(\beta;\cdot)$.
A {\em second generation I.F.S.} $\Phi$  consists of all maps of the form (\ref{eq-map1}), whose fixed points $\beta$ belong to the attractor $K_\Psi$ of the first generation I.F.S. $\Psi$:
\begin{equation}
\Phi = \{ \phi(\beta;\cdot), \; \beta \in K_\Psi \},
\label{eq-map1b}
\end{equation}

We use again eq. (\ref{eq-u1}) to define the operator $\bfuf$, replacing the set of maps $\Psi$ by $\Phi$. Since $K_\Psi$ is a compact set, the closure at r.h.s. of eq. (\ref{eq-u1}) is here redundant.
Let therefore $K_\Phi$ denote the fixed point of $\bfuf$:
$
K_\Phi=\bfuf\big(K_\Phi\big).
$
This set is the attractor of the {second generation I.F.S. $\Phi$ derived from $\Psi$ and $\alpha$}. We want to study its properties.

Our main result is the following: 
\begin{theorem}
For any finite, nonsingular hyperbolic I.F.S. $\Psi$
and for any $0<\alpha<1$ the attractor $K_\Phi$ of the second generation I.F.S. $\Phi$ derived from $\Psi$ and $\alpha$ consists of a finite union of closed intervals.
The same is true when $K_\Psi$ in definition (\ref{eq-map1b}) is replaced by a Cantor set $K$ admitting a construction of uniformly lower bounded dissection.
\label{teo-1}
\end{theorem}
Constructions of uniformly lower bounded dissection will be defined and described in the following.

The first part of this theorem has been conjectured in \cite{intj}, Conjecture 1, for disconnected, affine I.F.S. In the same work, a weaker result was found in a specific case: namely, it was proven (Theorem 1 \cite{intj}) that {\em when $\Psi$ is composed by a two--maps disconnected affine I.F.S. }(but the proof holds for any finite number of maps) {\em the attractor $K_\Phi$ contains an interval}.
Theorem \ref{teo-1} solves the problem completely (under the hypotheses above) and in wider generality.

It can now be used in conjunction with a localization analysis of the set $K_\Phi$. Formulae somehow simplify when the convex hull of $K_\Psi$ is the interval $[-1,1]$: $\cull(K_\Psi)=[-1,1]$. By a suitable rescaling we can always put ourselves in this situation. Then, it was proven in \cite{nalgo2,arxiv} (see also Lemma 1 in \cite{intj}) that $\cull(K_\Phi)=[-1,1]$ and
$ K_\Psi \subset K_\Phi \subset B_{2 \alpha}(K_\Psi)$, where $B_{2 \alpha}(K_\Phi)$ is the $2 \alpha$--neighborhood of $K_\Phi$.  Furthermore, consider the set $N_{\epsilon}$ introduced in \cite{intj}:
\begin{equation}
N_{\epsilon} = \{ x \mbox{ s.t. }
[x - \alpha - \epsilon, x + \epsilon] \cap (1-\alpha) K_\Psi = \emptyset \} .
\label{eq-neps}
\end{equation}
For any $\epsilon \geq 0$ $N_{\epsilon}$ is a finite collection of open intervals contained in the complement of $K_\Phi$. Theorem \ref{teo-1} proves that the set of gaps of $K_\Phi$, {\em i.e.}  $\cull(K_\Phi) \setminus K_\Phi$ is also a finite collection of open intervals, which includes $N_{\epsilon}$.

Finally, Theorem \ref{teo-1} proves that algorithm A2 in \cite{intj} terminates in a finite number of steps; hence it provides an efficient means of computation of the set $K_\Phi$.

From a more general perspective, the results of this paper belong to the study of the topological properties of sums of Cantor sets. A key ingredient of our proof is a result (Thm. \ref{lem-lem5} below) by Cabrelli {\em et al.} \cite{cabrelli} on finite sums of Cantor sets. In this context, the attractor $K_\Phi$ that we examine can be written as a geometric series of Cantor sets (see eq. (25) in \cite{intj}):
\begin{equation}
K_\Phi =  (1-\alpha) \sum_{j=0}^\infty  \alpha^j K_\Psi.
 \label{fan2}
\end{equation}

Not considered in the present paper, but to be mentioned in connection with these I.F.S., is the fact that {\em balanced measures} can also be associated to uncountable sets of maps \cite{elton,mendiv}. These measures have been studied in \cite{nalgo2,arxiv}, in the case of second generation I.F.S.: they are always of pure type and can be either absolutely continuous or singular continuous with respect to the Lebesgue measure on their support. Discriminating between the two cases appears to be an interesting problem.

This paper is organized as follows. In the next section we describe some basic properties of I.F.S. and their attractors. These properties are well known and we reproduce them here solely for convenience and as a way to introduce notations.

The successive Sect. \ref{sec-constr} introduces a standard way to describe Cantor sets in the real line. We mainly follow Sect. 2 in \cite{cabrelli} and we extend it by proving a few results needed in the remainder of the paper. The fundamental property of {\em uniform lower bounded dissection} ({\em ulbd}) of Cantor sets is also defined in this section. In Sect. \ref{sec-ulbd} we prove that this property holds for Cantor sets generated as attractors of two--maps non--singular hyperbolic I.F.S. We then derive some useful lemmas on the relation of ulbd property with certain operations on sets: section \ref{sec-union} contains an explicit construction by which it is proven that the union of two separated, ulbd Cantor sets is also ulbd, while Sect. \ref{sec-sum} proves that a ulbd Cantor set with prescribed properties can be found in the sum of a finite collection of ulbd Cantor sets.

We then move to the core of the problem: in Sect. \ref{sec-more} we recall Cabrelli {\em et al.} result on finite sums of ulbd Cantor sets, to which we add two consequent Lemmas.
This leads us to the final Sect. \ref{sec-proof} where we prove a proposition on finite sums of the geometric kind
(\ref{fan2}) and finally Theorem \ref{teo-1}.

\section{Basic properties of I.F.S. maps and attractors}
\label{sec-basic}

In this section we let $\Psi$ be a set of contractive transformations on $\bre$, not necessarily of the form (\ref{eq-contra}). When the cardinality of the set is finite, we will use the notation $\Psi:=\{\psi_i: i=1,\ldots,M\}$, and assume
there exist at least two maps with different fixed points. Then,
$m_{\Psi}=\min K_\Psi$ is strictly smaller than $M_{\Psi}=\max K_\Psi$. For all $n\ge 1$ let
$\Psi^n$ be the I.F.S. consisting of the $n$-fold composition of the maps in $\Psi$:

\begin{equation}
\Psi^n:=\big\{\psi_{i_1,\ldots,i_n}: i_1,\ldots,i_n=1,\ldots,M \big\},
\quad \psi_{i_1,\ldots,i_n}:=\psi_{i_1}\circ\cdots\circ \psi_{i_n}.
\label{eq-map1x}
\end{equation}

Then, clearly $\bfup^n=\bfu_{\Psi^n}\ \mbox{on}\ \ck$. Acting with $\bfup$ $n$ times on eq. (\ref{eqno(1)}) yields $K_\Psi = \bfup^n(K_\Psi) =\bfu_{\Psi^n} (K_\Psi)$
and this implies that,
for all $n\ge 1$
\begin{equation}
K_{\Psi^n}=K_\Psi.
\label{eqno(2)}
\end{equation}

\begin{lemma}
For a hyperbolic I.F.S., the following hold:

i) If $A\in\ck$ and $\bfup(A)\supseteq A$, then $K_\Psi\supseteq A$.

ii) If $B\in\ck$ and $\bfup(B)\subseteq B$, then $K_\Psi\subseteq B$.

iii) If $\Psi\subseteq\Psi'$ then $K_\Psi\subseteq K_{\Psi'}$.
\label{lem-1}
\end{lemma}

\smad i) The set $X:=\{D\in\ck: \bfup(D)\supseteq A\}$ is  a closed nonempty
(as $A\in X$) subset of $\ck$, thus a complete metric space with respect to the
Hausdorff metric. The map $\bfup$ is a contraction from $X$ into itself (since $D\in X$ implies $\bfup^2 (D) \supseteq \bfup (A) \supseteq A$), thus it has a
fixed point $C$: $\bfup(C)=C$. Because of uniqueness, this latter is the same as the I.F.S. attractor: $C=K_\Psi$. Moreover,
since $K_\Psi=C \in X$, then $K_\Psi=\bfup(K_\Psi)\supseteq A$.

ii) Same as i), with $X:=\{D\in\ck: \bfup(D)\subseteq B\}$.

iii) We have $\bfu_{\Psi'}\big(K_\Psi\big)\supseteq \bfup
\big(K_\Psi\big)=K_\Psi$. Now, iii) follows from i) with $A=K_\Psi$. \enpr

\begin{remark}
The previous lemma can be used to construct monotonic sequences of compact sets converging to the attractor $K_\Psi$: Take $A$ as in {\em i)} and define $A_n = \bfup^n (A)$, $B$ as in {\em ii)} and  $B_n = \bfup^n (B)$. Then,
\[
 A \subseteq \ldots \subseteq A_n \subseteq A_{n+1} \ldots \subseteq K_\Psi \subseteq \ldots \subseteq B_{n+1}
\subseteq B_{n} \subseteq \ldots \subseteq B.
\]
Typical choices for $A$ and $B$ are a finite set of fixed points, and $[m_{\Psi},M_{\Psi}]$, the convex hull of $K_\Psi$, respectively.
(Notice that an alternative proof of i) and ii) can be obtained via the previous equation.)
\end{remark}

A corollary of this result is the following
\begin{lemma}
\label{lem-2}
Let  $\Psi:=\{\psi_i: i=1,\ldots,M\}$. Let $\beta$ be the fixed point
of $\psi_{i_1,\ldots,i_n}$. Then $\beta\in K_\Psi$. Also $\psi_{i'_1,\ldots,i'_{k'}}(\beta)\in
K_\Psi$ for all $i'_1,\ldots,i'_{k'} =1,\ldots,M$.
\end{lemma}

\smad By definition, $\bfu_{\Psi^n}(\{\beta\})=\bfup^n(\{\beta\})\supseteq
\psi_{i_1,\ldots,i_n}(\{\beta\})=\{\beta\}$. Thus, by Lemma \ref{lem-1} i) $K_\Psi=K_{\Psi^n}\supseteq
\{\beta\}$. Next,

$$\psi_{i'_1,\ldots,i'_{k'}}(\beta)\in
\psi_{i'_1,\ldots,i'_{k'}}\big(K_\Psi\big)\subseteq
U_{\Psi^{k'}}\big(K_\Psi\big)=K_\Psi.    \ \enpr$$

Observe that the above also holds for infinite cardinality I.F.S.

\begin{lemma}
For finitely many maps there exists $\Psi'\subseteq\Psi$ such that $\Psi'$ has
precisely two elements and $m_{\Psi'}=m_{\Psi}$, $M_{\Psi'}=M_{\Psi}$.
\label{lem-3}
\end{lemma}

\smad Firstly, $m_{\Psi}\in K_\Psi=\bfup\big(K_\Psi\big)$. There exists
$\psi_1\in\Psi$ such that $$m_{\Psi}\in\psi_1\big( K_\Psi\big).$$ Thus, there exists
$x\in K_\Psi$ such that $m_{\Psi}=\psi_1(x)$. We have either $x=m_{\Psi}$ or
$x=M_{\Psi}$. In fact, in the opposite case there exist $x_1, x_2\in K_{\Psi}$ such
that $x_1<x<x_2$. Thus
  $\psi_1(x_1),\psi_1(x_2)\in \psi_1\big(K_{\Psi}\big)\subseteq
\bfup\big(K_{\Psi}\big)=K_{\Psi}$
  and one of the numbers   $\psi_1(x_1)$, $\psi_1(x_2)$ is less than
$m_{\Psi}$, a contradiction. Similarly, we can prove that there exist
$\psi_2\in\Psi$ and $y\in\{m_{\Psi}, M_{\Psi}\}$ such that $\psi_2(y)=M_{\Psi}$. Let
$\Psi':=\{\psi_1,\psi_2\}$. We claim that $\Psi'$ satisfies the Lemma. Note that,
since $K_{\Psi'}\subseteq K_\Psi \subseteq [m_{\Psi}, M_{\Psi}]$, it suffices to
prove that

\begin{equation}
m_{\Psi}, M_{\Psi}\in K_{\Psi'} \label{eq-eq(3)}
\end{equation}

We consider four different cases.

First case.  $x=m_{\Psi}$,  $y=M_{\Psi}$. Here, $m_{\Psi}$ is the fixed point of
$\psi_1$ and $M_{\Psi}$ is the fixed point of $\psi_2$. Now, eq. (\ref{eq-eq(3)}) follows from Lemma \ref{lem-2}.

Second case. $x=y=m_{\Psi}$. Then $m_{\Psi}$ is  the fixed point of $\psi_1$, and
$M_{\Psi}=\psi_2(m_{\Psi})$, and eq. (\ref{eq-eq(3)}) follows again from Lemma \ref{lem-2}.

Third case. $x=y=M_{\Psi}$. We proceed as in Second case.

Fourth case.  $x=M_{\Psi}$,  $y=m_{\Psi}$. We have
$\psi_1\circ\psi_2(m_{\Psi})=m_{\Psi}$, $\psi_2\circ\psi_1(M_{\Psi})=M_{\Psi}$, and eq. (\ref{eq-eq(3)})
follows again from Lemma \ref{lem-2}. \enpr

\section{Iterative construction of Cantor sets}
\label{sec-constr}

Recall that a Cantor set is a compact totally disconnected nonempty subset of $\bre$ with no
isolated points. Iterated function systems yield a convenient construction of families of Cantor sets on the real line. We now use a different description, of general scope, that has been also employed in \cite{cabrelli}. The main idea behind this construction is that the complement of a real Cantor set is a countable union of open intervals. How to organize this countable set is the core of the description, which requires symbolic coding, as follows.

Let $W$ be the set of finite binary words, with $\emptyset$ being the empty word:
$$W:={\emptyset}\cup\big\{w_1,\ldots,w_r: r=1,2,3,\ldots., w_i\in \{0,1\}\big\}.$$
Define the wordlength function $|\cdot|$ via $|\emptyset|=0$, $|w_1,\ldots,w_r|=r$. If $w,w'\in W$, let $ww'$ be the concatenation
of $w$ and $w'$, $w\emptyset=\emptyset w=w$.
Let us now associate a closed interval on the real line to each word in $W$: that is to say, we define a map $\wii : W \rightarrow \{ \mbox{closed intervals } \subset \bre \}$, such that
\begin{equation}
\wii : w \rightarrow  I_w := [a_w,b_w],
\label{eq-defwii}
\end{equation}
with $a_w$ and $b_w$ denoting the end points of $I_w$ (clearly, we always require $a_w < b_w$). This map is defined iteratively. The initial seed is
an arbitrary interval $I=[a,b]$ that is associated to the empty word. That is, $I_{\emptyset}=I$, $a_{\emptyset}=a$, $b_{\emptyset}=b$. The iteration rule is then the following: given $I_w=[a_w,b_w]$ with $a_w<b_w$, choose two points $c_w$ and $d_w$, so that $a_w<c_w<d_w<b_w$ holds with strict inequalities and define the new intervals $I_{w0}=[a_{w0},b_{w0}]$ and  $I_{w1}=[a_{w1},b_{w1}]$ via
$$a_{w0}=a_w, b_{w0}=c_w, a_{w1}=d_w, b_{w1}=b_w.$$
In simpler terms, the interval corresponding to a word $w$, of length $|w|$ generates two intervals, corresponding to the words $w0$ and $w1$ of length $|w|+1$. It is  convenient to define the {\em ratios of dissection} $r(\wii)_{w0}$ and $r(\wii)_{w1}$ associated to these intervals, as
\begin{equation}
r(\wii)_{wj}=\disp{b_{wj}-a_{wj}\over b_w-a_w}={d(I_{wj})\over d(I_w)}, \;
j = 0,1
\label{eq-crat}
\end{equation}
where here and in the following we use the notation
$d(A)=\mbox{diam }(A)$ for the diameter of the set $A$. Clearly, $d(A) = \max A-\min A$ for every nonempty compact subset $A$ of $\bre$. Note that, since $d(I_{w0})+d(I_{w1})<d(I_{w})$, we have
\begin{equation}
\label{eq-eq5}
r(\wii)_{w0}+r(\wii)_{w1}<1.
\end{equation}
The ratio of dissection $r(\wii)_{v}$ is so defined for any word $v\in
W\setminus\{\emptyset\}$ and it measures the ratio of the diameters of $I_v$ and of its immediate ancestor $I_{v'}$ (associated with the word $v'$, obtained from $v$ by deleting the last binary digit).

To complete the construction, take the union of the intervals $I_w$ of fixed length $|w|$, and then intersect these latter sets:
\begin{equation}
\label{eq-eq4x}
C_n(\wii):=\bigcup\limits_{|w|=n} I_w, \quad C(\wii):=\bigcap\limits_{n=0}^{\infty}
C_n(\wii).
\end{equation}
Note that $I_{wj}\subseteq I_w$, thus in particular
$C_{r+1}(\wii)\subseteq C_r(\wii)$ and $C(\wii)$ is not empty.
We say that $C(\wii)$  {\it is the quasi-Cantor set constructed on $\wii$} or that $\wii$ {\em constructs} $C(\wii)$.

Consider the following specific case: For $j=0,1$, let $j^n$ denote the word composed of the letter $j$ repeated $n$ times,
with $n \in \bna$, where $j^0=\emptyset$. These words are labels of the extreme intervals in $C_n(\wii)$, so that
$$\min I_{0^n}=\min I, \quad \max I_{1^n}=\max I,$$
and
$$\min C(\wii)=\min I, \ \ \max C(\wii)=\max I.$$

The correspondence $\wii \rightarrow C(\wii)$ is not one-to-one: the same Cantor set may be constructed on different $\wii$'s. This is clearly seen by considering  the complement of a quasi-Cantor set $C(\wii)$: define the {\it gaps} of $C(\wii)$ as the bounded connected components of the complement of $C(\wii)$. Their countable union is precisely $I_\emptyset \setminus C(\wii)$.

\begin{lemma}
\label{lem-gaps}
The gaps of $C(\wii)$ are the sets
\begin{equation}
\label{eq-eq6}
I_w\setminus \Big(I_{w0}\cup I_{w1}\Big): w\in W.
\end{equation}
\end{lemma}

\smad The sets  in (6) are clearly bounded components of the complement of $C(\wii)$.
Conversely, suppose $A$ is a bounded component of the complement of $C(\wii)$ and take
$x\in A$. Then $x\in I_\emptyset$. In fact, in the opposite case, either $x<a$ or $x>b$. In
the former case $x\in ]-\infty,a[$, in the latter $x\in]b,+\infty[$, so that $A$
being the connected component of $x$ contains either $]-\infty,a[$ or
$]b,+\infty[$, thus is unbounded, a contradiction. Now recall that $I_\emptyset = C_0(\wii)$, so that $x\in
C_0(\wii)\setminus\Big(\bigcap\limits_{r=0}^{+\infty} C_r(\wii)\Big)$ and there exists $r\in\bna$ such that $x\in C_r(\wii)\setminus C_{r+1}(\wii)$, and also $w\in W$ such that $|w|=r$ and $x\in I_w\setminus
\Big(I_{w0}\cup I_{w1}\Big)$. Therefore, $A=I_w\setminus \Big(I_{w0}\cup
I_{w1}\Big)$. \enpr

Note that $I_w\setminus \Big(I_{w0}\cup I_{w1}\Big)=]c_w,d_w[$. Therefore, the above construction of $C(\wii)$, defined by  $\wii$, can also be seen as a construction of its complementary in $[a,b]$, described by a function ${\cal G}$ from $W$ to the set of open intervals. Keeping fixed the image of this map, {\em i.e.} the gaps, any map ${\cal G}$, that respects a simple prescription (gaps appear in interlacing sequence) yields the same Cantor set.
We will use this freedom later on in the paper. We will also use a specific symbol for the diameter of the gaps:
\begin{equation}
\label{eq-gapsize}
\gamma(\wii)_w:=d_w-c_w=\min I_{w1}-\max I_{w0},
   \quad\gamma(C)=\sup\limits_{w\in W} \gamma(\wii)_w.
   \end{equation}

We now give a condition for $C(\wii)$ being a Cantor set. Let us start with a symbolic coding of all points in $C(\wii)$. Denote by $\witi W$ the set of
infinite strings of $0$ and $1$, {\em i.e.} $\witi W = \{0,1\}^{\bna \setminus \{0\}}$, and for $\tilde w\in \witi W$, write $\tilde w= i_1 i_2 i_3 \cdots$. Also, let $\tilde w_n$ be the finite string of length $n$ obtained by truncation of $\tilde w$: $\tilde w_n=i_1 \cdots i_n$. With this vocabulary, a point $x \in \bre$ belongs to $C(\wii)$ if and only if there exists $\tilde w\in \witi W$ such that
$x\in\bigcap\limits_{n=1}^{\infty} I_{\tilde w_n}$. Indeed, when the set is Cantor, this intersection is the singleton $\{x\}$, as the following
standard lemma shows:

\begin{lemma}
\label{lem-cantor}
The set $C(\wii)$ is a  Cantor set if
and only if for every $\tilde w\in\witi W$
\[ d(I_{\tilde
w_n}) \llim\limits_{n\to +\infty} 0.
\]
\end{lemma}

\smad
In fact, being the sets $I_{\tilde w_n}$ serially enclosed, the sequence of their diameters is monotonic, and if for a certain $\tilde w\in\witi W$
it does not converge to zero, then it has a strictly positive limit: $d(I_{\tilde w_n})\llim\limits_{n\to +\infty} c>0$. In this case, $C(\wii)\supseteq
\bigcap\limits_{n=1}^{\infty} I_{\tilde w_n}\supseteq [\alpha,\beta]$ with $\alpha<\beta$. Thus $C(\wii)$
is not totally disconnected.

Conversely, suppose that for every $\tilde w\in\witi W$ we have $d(I_{\tilde w_n})\llim\limits_{n\to +\infty} 0$. Let $x\in C(\wii)$ and take $\tilde w\in \witi W$
such that $x\in\bigcap\limits_{n=1}^{\infty} I_{\tilde w_n}$. Put $I_{\tilde
w_n}=[a_n,b_n]$. Since $a_n\le x\le b_n$ and $b_n-a_n\llim\limits_{n\to +\infty}0$,
for every $U$ neighborhood of $x$ there exists $n$ such that $[a_n,b_n]\subseteq
U$. If $\tilde w'\in \witi W$ and $\tilde w'_n=\tilde w_n$ but $\tilde w'\ne \tilde
w$ and $x'\in\bigcap\limits_{n=1}^{\infty} I_{\tilde w'_n}$, then $x'\in C(\wii)\cap
[a_n,b_n]$,   thus $x'\in U$, but $x'\ne x$, so that $C(\wii)$ has no isolated points.
Also, since $C_n(\wii)$ is the {\em disjoint} union of closed intervals including $[a_n,b_n]$,
then the component of $x$ in $C(\wii)$ is contained in $[a_n,b_n]$, for every $n$, thus
in $\bigcap\limits_{n=1}^{\infty} I_{\tilde w_n}$, which, since
$b_n-a_n\llim\limits_{n\to + \infty} 0$ amounts to $\{x\}$, hence $C(\wii)$ is totally disconnected.~\enpr

We can use the previous lemma in conjunction with the following:
\begin{lemma}
\label{lem-cant}
A sufficient condition for $C(\wii)$ being a Cantor set is that there
exists a positive constant $a$, such that  all dissection ratios are larger than, or equal to $a$.
\end{lemma}

\smad
Equation (\ref{eq-eq5}) implies that
$r(\wii)_{wj}<1 -r(\wii)_{w(1-j)}\le 1-a<1$.
Hence, since $d(I_{wj})=r(\wii)_{wj} d(I_{w})\le (1-a)d(I_{w})$, it follows that
for every $\tilde w \in \tilde W$
$d(I_{\tilde w_n})\le (1-a)^n d(I)$, thus  $d(I_{\tilde w_n})\llim\limits _{n\to +\infty} 0$. \enpr

Clearly, this condition is not so much intended to exclude that  dissection ratios get too small (which could still be compatible with having a Cantor set, and hence the condition is not necessary), rather, because of the inequality (\ref{eq-eq5}) it implies that dissection ratios cannot tend to one.

Given the role that this condition will play in the following, we find it convenient to embody it into a formal definition:
\begin{definition}
\label{def-lowbd}
A construction $\wii$ that satisfies the condition in Lemma \ref{lem-cant}
will be said to be of {\em uniformly lower bounded dissections (ulbd)}, and a Cantor set admitting one such construction will also be said to possess the same property.
\end{definition}

We end this section by defining a further element in the algebra of quasi-Cantor sets. Observe that each interval in the above construction can be thought of as the starting interval in the construction of a Cantor set, subset of the former. In fact, let the mapping $\wii$ be fixed, and let us focus on a finite word $w \in W$ and on its associated interval $I_w=\wii(w)$. Define a new mapping $\wii_w$ by the formula (compare with eq. (\ref{eq-defwii})):
\begin{equation}
\wii_w : w' \rightarrow \wii_w(w') := \wii( w w') = I_{w w'}.
\label{eq-defwii2}
\end{equation}
Denote by $C(\wii_w)$
the quasi Cantor set constructed in the set $I_w$ by this mapping.
Clearly
\begin{equation}
C({\cal I}_w)=I_w\cap C({\cal I}).
\label{eq-ciccia}
\end{equation}
Moreover, the set of gaps of $\wii_w$ is contained in
the set of gaps of $\wii$, and the set of ratios of $\wii_w$ is contained
in the set of ratios of $\wii$. Namely, $\gamma
(\wii_w)_{w'}=\gamma(\wii)_{ww'}$
  and $r(\wii_w)_{w'}=r(\wii)_{ww'}$. Of course, it also holds true that

  $$\max C(\wii_w)=\max I_w, \quad \min C(\wii_w)=\min I_w.$$

\section{Ulbd property of I.F.S. Cantor sets}
\label{sec-ulbd}
In this section we prove that Cantor sets constructed via I.F.S. with two maps of the form (\ref{eq-contra}) possess the uniform lower bounded dissection property. Prior to that, we need a technical lemma.
\begin{lemma}
\label{lem-issimo}
Let $\Psi$ be a set of $C^2$ maps from an interval $I$ to itself, as in eq. (\ref{eq-map00}) that satisfy condition (\ref{eq-contra}).
For any finite word $w \in W$, define $\psi_w = \psi_{w_1} \circ \cdots \circ \psi_{w_n}$, as in eq. (\ref{eq-map1x}).
Then there exists $c>0$ such that, for any interval $J \subseteq I$, any $w \in W$, any $i = 1,\ldots,M$, we have
\begin{equation}
 r_{wi} :=  \frac{d(\psi_{wi}(J))}{  d(\psi_{w}(J))}\ge c.
\label{eq-crat2}
\end{equation}
\end{lemma}

\smad
For any $w \in W$ and any $i$ we need to estimate the ratios $r_{wi}$.
Let us first consider the numerator in eq. (\ref{eq-crat2}): this is the length of an interval, that can be evaluated as
\[
d(\psi_{wi}(J)) = |\psi_{wi}'(\eta)| d(J),
\]
where $\eta$ is a point in $J$. Similarly,
$
d(\psi_{w}(J)) = |\psi_{w}'(\zeta)| d(J),
$
$\zeta \in J$. The chain rule for the derivative of these composed functions leads us to define two sequences of points $\eta_k$, $\zeta_k$, for $k=1,\ldots,n$, as follows:
\[
\eta_k = (\psi_{w_{k+1}} \circ \cdots \circ \psi_{w_n} \circ \psi_i) (\eta), \;\; k = 1,\ldots,n-1 ,
\]
\[
\zeta_k = (\psi_{w_{k+1}} \circ \cdots \circ \psi_{w_n}) (\zeta), \;\; k = 1,\ldots,n-1 ,
\]
and $\eta_n = \psi_j(\eta)$, $\zeta_n=\zeta$.
With these notations, the derivative of the composed function can be factored as
\[
\psi_{wi}'(\eta) = \psi'_{w_1} (\eta_1) \cdot \ldots \cdot \psi'_{w_n} (\eta_n) \cdot \psi_i'(\eta),
\]
\[
\psi_{w}'(\zeta) = \psi'_{w_1} (\zeta_1) \cdot \ldots \cdot \psi'_{w_n} (\zeta_n) .
\]
Because of contractivity of the maps, the points $\eta_k$, $\zeta_k$ approach each other geometrically, when $n$ grows. In fact, $\eta_n, \zeta_n \in J$ and
\[
  \eta_k, \zeta_k \in (\psi_{w_{k+1}} \circ \cdots \circ \psi_{w_n}) (J), \; k = 1,\ldots,n-1.
\]
Using eq. (\ref{eq-contra}) we obtain
\[
  d((\psi_{w_{k+1}} \circ \cdots \circ \psi_{w_n}) (J)) \leq \delta^{n-k} d(J), \; k = 1,\ldots,n-1.
\]
The above information permits to compute the logarithm of the inverse of the ratio $r_{wi}$: we call it $l_{wi}$ and we show that it is bounded from above. In fact,
\begin{equation}
l_{wi}:=-\log(r_{wj}) =
- \log(|\psi_i'(\eta)|) + \sum_{k=1}^n \log (|\psi'_{w_k} (\zeta_k)|) - \log (|\psi'_{w_k} (\eta_k)|).
\label{eq-crat3}
\end{equation}
Therefore,
\begin{equation}
l_{wi} \leq \log(1/\sigma)
+ \sum_{k=1}^n |\log (|\psi'_{w_k} (\zeta_k)|) - \log (|\psi'_{w_k} (\eta_k)|)|.
\label{eq-crat4}
\end{equation}
Consider now the functions $g_i(x) = \log(|\psi'_i(x)|)$, where $i=0,1$.
Because of  eq. (\ref{eq-contra}) $g_i(x)$ is differentiable and
\[
  |g'_i(x)| = \frac{|\psi''_i(x)|}{|\psi'_i(x)|},
\]
so that each term in the summation at r.h.s. of eq. (\ref{eq-crat4}) can be estimated as
\begin{equation}
 |\log (|\psi'_{w_k} (\zeta_k)|) - \log (|\psi'_{w_k} (\eta_k)|)| =  \frac{|\psi''_{w_k}(\theta_k)|}{|\psi'_{w_k}(\theta_k)|} |\zeta_k-\eta_k| \leq
 \frac{B}{\sigma} \delta^{n-k} d(J),
\label{eq-crat5}
\end{equation}
with $\theta_k$ an intermediate point between $\eta_k$ and $\zeta_k$ and where $B$ is the maximum of the absolute value of the second derivative of all $\psi_i$'s over $I$. In conclusion, we have
\begin{equation}
l_{wi} \leq \log(1/\sigma)
+ \frac{B}{\sigma} d(J) \sum_{k=1}^n  \delta^{n-k} \leq
\log(1/\sigma) + \frac{B}{\sigma(1-\delta)} d(I).
\label{eq-crat6}
\end{equation}
The term at r.h.s. is a finite quantity $C$, independent of $w$ and $i$, and this proves the lemma: $r_{wi} \geq e^{-C}$ for all $w \in W$, $i=0,\ldots,M$. \enpr

The main lemma of this section is now the following.
\begin{lemma}
\label{lem-nuovo}
Let $\Psi = \{\psi_0,\psi_1\}$ be a set of $C^2$ I.F.S. maps that satisfy condition (\ref{eq-contra}).  Let these maps have different fixed points.
Then $K_\Psi$ is either an interval, or a Cantor set that admits a construction $\wii$ that is of ulbd.
\end{lemma}

\smad
Remark that the attractor of a two--maps I.F.S. with different fixed points is either a full interval or a Cantor set, as it is easy to see. In fact, let $I_\emptyset = \cull(K_\Psi) = [a,b]$, the convex hull of $K_\Psi$.
Let now $J_i = \psi_i([a,b])$, $i=0,1$. The extreme point $a$ must belong to one of these two intervals, and $b$ to the other: in fact, they belong to $K_\Psi$ and therefore also to $U_{\Psi}(K_{\Psi})$. If these two intervals are not disjoint, we have that $U_{\Psi}([a,b]) = [a,b]$ and therefore the attractor is the full convex hull, $K_\Psi = [a,b]$.
In the opposite case, $J_0$ is either strictly to the left of $J_1$, or to its right. In the first case we assign a permutation $g$ of $\{0,1\}$, such that $g(0)=0$ and $g(1)=1$ (the identity). In the second case we invert indices: $g(0)=1$, $g(1)=0$, so that in both cases we define $I_i = [a_i,b_i] = \psi_{g(i)}(I_\emptyset)$ and we have $a_0=a$, $b_1=b$. Disconnectedness of the two intervals imply that $b_0 < a_1$, thereby completing the first step in the construction of the Cantor set.

We then proceed by induction: consider the words $w \in W$ of length $n-1$, 
the maps $\psi_w$ (as in eq. (\ref{eq-map1x})) and the permutation $g$ of the set of $n-1$ letter words that defines the lexicographically ordered intervals $I_w = [a_w,b_w] = \psi_{g(w)}(I_\emptyset)$.
Define the intervals $J_{wi} = (\psi_{g(w)} \circ \psi_i)(I_\emptyset)$, for $i=0,1$. Clearly, $J_{wi} \subset I_w$, and these two intervals are disjoint. Extend the permutation $g$ to the set of $n$-letter words as follows: $g(w0) = g(w) 0$, $g(w1) = g(w) 1$ if $J_{w0}$ is to the left of $J_{w1}$, or $g(w0) = g(w) 1$, $g(w1) = g(w) 0$ if otherwise. This implies that $I_{wi}=[a_{wi},b_{wi}] = \psi_{g(wi)}(I_\emptyset)$, and $a_w=a_{w0}<b_{w0}<a_{w1}<b_{w1}=b_w$.
This proves that the map $\wii$ so defined yields a Cantor set.

Let us now prove that this construction is of {\em ulbd}. In fact, 
contraction ratios are defined by eq. (\ref{eq-crat}): in this case, they are given by
\begin{equation}
r(\wii)_{wj}=\frac{d(I_{wj})}{d(I_w)} =
 \frac{d(\psi_{g(wj)}(I_\emptyset))}{  d(\psi_{g(w)}(I_\emptyset))},
\label{eq-crat2x}
\end{equation}
with $j=0,1$. Since $g(wj) = g(w) h_w(j)$, where $h_w$ is a permutation of a last letter (that depends on $w$, but this is not an issue), we can apply Lemma \ref{lem-issimo} to prove that these ratios are uniformly bounded from below.
\enpr

\section{Union of Cantor sets}
\label{sec-union}

This section explains how to organize the union of two ulbd Cantor sets, into a single construction, $\wii'$, that is also ulbd, perhaps with a smaller lower bound. To do this, we shall exploit the non uniqueness of the construction.

\begin{lemma}
\label{lem-lem7}
Let $C^{(1)}$ and $C^{(2)}$ be two  Cantor sets of ulbd, separated so that $\max C^{(1)}<\min C^{(2)}$. Then, $C:=C^{(1)}\cup C^{(2)}$ is a Cantor set, that admits a construction $\wii$ that is also of ulbd. Moreover,
\begin{equation}
\label{eq-eq11}
\gamma(C) =
  \max\big\{\gamma(C^{(1)}), \gamma(C^{(2)}), \min C^{(2)}-\max
C^{(1)}\big\}.
\end{equation}

Finally, lower bounds can be estimated as follows. Let $a^{(1)}$ and $a^{(2)}$ be the (strictly positive) lower bounds to the dissection rates of $C^{(1)}$ and $C^{(2)}$. Let also
\begin{equation}
\label{eq-eq10}
a = \min\Big\{\disp{ {\max C^{(1)}-\min C^{(1)}\over \min C^{(2)}-\min
C^{(1)}}}, \disp{ {\max C^{(2)}-\min C^{(2)}\over \max C^{(2)}-\max C^{(1)}}},a^{(1)},a^{(2)}
\Big\}.
\end{equation}
Then, there exists a universal value $a'\in]0,a]$, that depends only on $a$ and not on the Cantor sets $C^{(1)}$ and $C^{(2)}$, so that $C$ admits a ulbd construction, with dissection ratios larger than, or equal to $a'$, with
\begin{equation}
\label{eq-aprime}
a'=\min\big\{ {a\over 2}, {a^2\over a+1} \big\}.
\end{equation}
\end{lemma}
  \smad
First, it is clear that $C:=C^{(1)}\cup C^{(2)}$ is a Cantor set, because of the separation condition $\max C^{(1)}<\min C^{(2)}$. Therefore, it can be constructed on a map $\wii$, although not in a unique way. Since the set of gaps do not depend on the construction $\wii$, the gaps of $C$ are the union of those of $C^{(1)}$ and $C^{(2)}$ and the open interval $]\max C^{(1)},\min C^{(2)}[$, eq. (\ref{eq-eq11}) follows. We now need to prove that such a construction exists, that is of uniformly lower bounded dissection. We denote  by $\wii$ this construction.

Suppose $C^{(1)}$ and $C^{(2)}$ are constructed on $\wii^{(1)}$ and $\wii^{(2)}$, with ratios of dissection $r(\wii^{(1)})_w$,  $r(\wii^{(2)})_w$, which by hypothesis are all larger than, or equal to $a$.
Without loss of generality, assume that $I^{(1)}$ is wider than $I^{(2)}$:
\begin{equation}
\label{eq-eq12}
d(I^{(1)})\ge d(I^{(2)}).
\end{equation}
The hypothesis and eq. (\ref{eq-eq5}) imply that $r(\wii^{(j)})_w< 1-a$ for every $w\in W\setminus\{\emptyset\}$ and $j=1,2$.
Then,
$d(I^{(j)} _{wi})=r^{(j)}_{wi}\, h (I^{(j)} _w)\le (1-a) h
(I^{(j)}_w)$ and we conclude that
$$d(I^{(j)}_w)\le (1-a)^{|w|} d(I^{(j)})$$ for each $w\in W$ and
$i,j=1,2$. Then, because of (\ref{eq-eq12}), there exists $\bar n \in \bna$ such that
\begin{equation}
\label{eq-eq13}
 d(I^{(1)}_{1^n})\ge d(I^{(2)}) \ \ \ \mbox{for}\ n\le \bar n,\qquad
   d(I^{(1)}_{1^n})< d(I^{(2)}) \ \ \ \mbox{for}\ n> \bar n.
\end{equation}
In view of this result, let us define a construction $\wii$ as follows:
for any $w \in W$ define $I_w$ by
\begin{equation}
\label{eq-constr}
I_w= \left\{ \begin{array}{ll}
   \big[\min I^{(1)}_w, \max I^{(2)}\big] & \mbox{if}\ w=1^n,  0 \le n\le\bar
n
   \quad\mbox{(first\ case)}\\
  I^{(1)}_{1^{\bar n}\, w'} &\ \mbox{if}\ w=1^{\bar n}0w'
  \quad\quad \quad\ \mbox{(second\ case)}\\
  I^{(2)}_{w'} &\ \mbox{if}\ w=1^{\bar n}1w'
  \quad\quad \quad\ \mbox{(third\ case)}\\

   I^{(1)}_w&\ \mbox{otherwise.}
   \quad\quad\quad \quad\ \mbox{(fourth\ case)}\end{array} \right.
  \end{equation}

We first prove that $C^{(1)}\cup C^{(2)}$ is constructed on $\wii$. That is,
$\wii$ constructs a quasi Cantor set and, putting $C_n:=
\bigcup\limits_{|w|=n} I_w$, we have
\begin{equation}
\label{eq-eq14}
C^{(1)}\cup C^{(2)}=\bigcap\limits_{n=0}^{\infty} C_n.
\end{equation}
We will use systematically the following evident remark:
Since $\max I^{(1)}=\max C^{(1)}<\min C^{(2)}=\min I^{(2)}$
and the intervals $I^{(j)}_w$ are all contained in $I^{(j)}$ ($j=1,2$),
then any element of $I^{(1)}_w$ is strictly less than any element of $I^{(2)}_{w'}$ for every $w,w'\in W$.

Note that by hypothesis $\min I^{(1)}<\max I^{(1)} <\min I^{(2)}<\max I^{(2)}$ and in our construction $\wii$, eq. (\ref{eq-constr}), $I=I_{\emptyset}=[\min I^{(1)}, \max  I^{(2)}]$. Let again $I_w=[a_w,b_w]$. Note that by
construction, in any case $a_w<b_w$. If $w=1^n, \ \ n\le\bar n$, this follows from the above remark, since $\min I^{(1)}_w  \le \max I^{(1)}$. In the other cases the intervals
considered are of the form $I^{(j)}_v$ with $j=1,2$, which by hypothesis satisfy
the inequality $a_v<b_v$.

Next, we have to prove that $a_w=a_{w0}<b_{w0}<a_{w1}<b_w=b_{w1}$, {\em i.e.}
\begin{eqnarray}
 \min I_{w0} & = \min I_{w},
\label{eq-eq15}
\\
 \max I_{w0}& <\min I_{w1},
\label{eq-eq15'}
\\
\max I_{w1} & =\max I_{w}.
\label{eq-eq15''}
\end{eqnarray}

In first case of eq. (\ref{eq-constr}) we have $
I_w=   \big[\min I^{(1)}_w, \max I^{(2)}\big]$
and either one of the two possibilities holds:
\begin{eqnarray}
 \label{eq-eq16'}
 n<\bar n, &  \quad
I_{w0}=I^{(1)}_{w0}, & \quad I_{w1}=\big[\min I^{(1)}_{w1}, \max I^{(2)}\big];
\\
  n=\bar n, & \quad I_{w0}=I^{(1)}_{1^{\bar n}}=
I^{(1)}_w, &  I_{w1}= I^{(2)}.
\label{eq-eq16''}
\end{eqnarray}
We see that (\ref{eq-eq15}) holds in both cases.
Moreover, (\ref{eq-eq15'}) holds trivially if $n<\bar n$, while if $n=\bar n$ we have
$\max I_{w0}=\max I^{(1)}_w\le \max I^{(1)}<\min I^{(2)}= \min
I_{w1}.$

Finally (\ref{eq-eq15''}) is trivial in both subcases. In second case,  $w0=1^{\bar n}0w'0$ and $w1=1^{\bar n}0w'1$. Thus,
\begin{equation}
\label{eq-eq17}
I_w=I^{(1)}_{1^{\bar n}\, w'}, \quad
  I_{w0}=I^{(1)}_{1^{\bar n}\, w'0},\quad I_{w1}=I^{(1)}_{1^{\bar n}\,
w'1},
\end{equation}

so that (\ref{eq-eq15}), (\ref{eq-eq15'}) and (\ref{eq-eq15''}) follow from the corresponding properties of $I^{(1)}_w$. In the third case, then $w0=1^{\bar n}1w'0$ and $w1=1^{\bar n}1w'1$.
Thus,
\begin{equation}
\label{eq-eq18}
I_w= I^{(2)}_{w'}, \quad  I_{w0}= I^{(2)}_{w'0}, \quad  I_{w1}=
I^{(2)}_{w'1}, \end{equation}

and we proceed as above. In the fourth case we have
\begin{equation}
\label{eq-eq19}
I_w=I^{(1)}_w,\quad I_{w0}=I^{(1)}_{w0}, \quad I_{w1}=I^{(1)}_{w1},
\end{equation}
and we proceed similarly.
  Thus, we have proven that in fact $\wii$ constructs a quasi Cantor set. We now
prove that $C^{(1)}\cup C^{(2)}$ is constructed on $\wii$, that is eq. (\ref{eq-eq14}). Note that for every $n>\bar n$ we have
\begin{equation}
\label{eq-eq20}
C^{(1)}_{n}\cup C^{(2)}_{n} \subseteq
C_n\subseteq C^{(1)}_{n-\bar n-1}\cup C^{(2)}_{n-\bar n-1}.
\end{equation}

To prove (\ref{eq-eq20}), suppose first $x\in C^{(1)}_{n}\cup C^{(2)}_{n}$. Then, either $x\in
C^{(1)}_{n}$ or $x\in  C^{(2)}_{n}$. In the former case, there exists $w\in W$ with
$|w|=n$ such that $x\in I^{(1)}_w$, and either $w$ can be written as $w=1^{\bar n}w'$,
in which case $x\in I_{1^{\bar n}0w'} \subseteq C_{n+1}\subseteq C_n$, or $w$ is
not of the form $w=1^{\bar n}w'$, in which case $x\in I_w\subseteq C_n$. If instead
$x\in  C^{(2)}_{n}$, then $x\in I^{(2)}_w$ for some $w\in W$ with $|w|=n$. Then,
$x\in I_{1^{\bar n}1w}\subseteq C_{n+\bar n+1}\subseteq C_n$. The first inclusion
in (\ref{eq-eq20}) is so proven. Let us prove the second. Suppose that $x\in C_n$, thus $x\in I_w$ for
some  $w\in W$ with $|w|=n$. Then by definition, since $n>\bar n$ we are not in the
first case in definition of $I_w$. If the second case holds, then $x\in
I^{(1)}_v$ with $|v|=n-1$; in the third case one has $x\in I^{(2)}_v$ with $|v|=n-1-\bar n$, and in the fourth case $x\in I^{(1)}_v$ with $|v|=n$. In any case, since the sequences of sets
$C^{(j)}_n$ are decreasing we have that $x\in C^{(1)}_{n-\bar n-1}\cup C^{(2)}_{n-\bar n-1}$ and (\ref{eq-eq20}) is proven.

At this point, from (\ref{eq-eq20}), since $C^{(j)}_n\subseteq C^{(j)}_0=I^{(j)}$
and $I^{(1)}\cap I^{(2)}=\emptyset$, and $\bigcap\limits_{n=0}^{\infty} C^{(j)}_n=
C^{(j)}$, eq. (\ref{eq-eq14}) easily follows. Thus, we have proven that $C^{(1)}\cup C^{(2)}$ is constructed on $\wii$.

Finally, we have to prove the fundamental part of the lemma, that is, there exists a positive constant $a'$ such that $r(\wii)_{wj}\ge a'$ for all $w\in W$, $j=0,1$. Note that by the hypothesis (\ref{eq-eq10}) we have

$$\max I_2-\min I_2\ge a\big(\max I_2-\min I_2+\min I_2-\max I_1\big)$$

hence, using also (\ref{eq-eq13}),
\begin{equation}
\label{eq-eq21}
\max I^{(1)}_{1^n}-\min I^{(1)}_{1^n}\ge
\max I_2-\min I_2\ge {a\over 1-a} \big(\min I_2-\max I_1\big)
\end{equation}

for every $n\le \bar n$.  Similarly, since eq. (\ref{eq-eq10}) implies that

$$\max I^{(1)}-\min I^{(1)}\ge a\big(\min I^{(2)}-\max I^{(1)}+\max
I^{(1)}-\min I^{(1)}\big)$$

we have
\begin{equation}
\label{eq-eq21'}
\max I^{(1)}-\min I^{(1)}\ge {a\over 1-a} \big(\min I^{(2)}-\max
I^{(1)}\big).
\end{equation}

Moreover, by (\ref{eq-eq13}),
\[
d(I^{(2)}) > d(I^{(1)}_{1^{\bar n+1}}) = r(\wii^{(1)})_{1^{\bar n+1}}\,
d(I^{(1)}_{1^{\bar n}}) \geq a d(I^{(1)}_{1^{\bar n}}) .
\]

Hence,
\begin{equation}
\label{eq-eq22}
d(I^{(2)})\le d(I^{(1)}_{1^{\bar n}})\le {1\over a} d(I^{(2)}).
\end{equation}
Following these inequalities, we can evaluate $r(\wii)_{wj}$, with $w\in W$.

When (\ref{eq-eq16'}) holds, we first estimate $r(\wii)_{w0}$:
\begin{equation}
\label{eq-eq23}
r(\wii)_{w0}={\max I_{w0}-\min I_{w0}\over  \max I_{w}-\min I_{w}}=
{\max I^{(1)}_{w0}-\min I^{(1)}_{w0}\over  \max I^{(2)}- \min I^{(1)}_{w}}. \end{equation}

Now, since $w=1^n$, $n\le \bar n$, by (\ref{eq-eq21}) we have

\begin{eqnarray}
\max I^{(2)}- \min I^{(1)}_{w}  = &
\max I^{(2)}-\min I^{(2)}
+\min I^{(2)}-\max  I^{(1)}_{1^n}+ \max I^{(1)}_{1^n}-\min I^{(1)}_{1^n}
\nonumber \\
= & \max I^{(2)}-\min I^{(2)}+\min I^{(2)}-\max  I^{(1)}+
  \max I^{(1)}_{1^n}-\min I^{(1)}_{1^n}
\nonumber \\
\le & 2\big(\max I^{(1)}_{1^n}-\min I^{(1)}_{1^n}\big)+\min I^{(2)}-\max
I^{(1)}
\nonumber \\
\le & \Big(2+{1-a\over a}\Big) \big(\max I^{(1)}_{1^n}-\min
I^{(1)}_{1^n}\big)
\nonumber \\
= &\Big(1+{1\over a}\Big) \big(\max I^{(1)}_{w}-\min I^{(1)}_{w}\big)
={1\over r(\wii^{(1)})_{w0}}   {a+1\over a}
  \big(\max I^{(1)}_{w0}-\min I^{(1)}_{w0}\big)
\nonumber \\
\le & {a+1\over a^2} \big(\max I^{(1)}_{w0}-\min I^{(1)}_{w0}\big).
\label{eq-manylin}
\end{eqnarray}

Thus, also in view of (\ref{eq-eq23}),  we have
\begin{equation}
\label{eq-eq24}
r(\wii)_{w0}\ge {a^2\over a+1}.
\end{equation}

Let us now consider $r(\wii)_{w1}$, still when (\ref{eq-eq16'}) holds. We have

$$r(\wii)_{w1}={\max I_{w1}-\min I_{w1}\over \max I_{w}-\min I_{w}}=
{\max I^{(2)}-\min  I^{(1)}_{w1}\over \max I^{(2)}-\min  I^{(1)}_{w} }.$$

Now,

$$\max I^{(1)}_{w}-\min I^{(1)}_{w1}=\max I^{(1)}_{w1}-\min
I^{(1)}_{w1}=$$
$$r(\wii^{(1)})_{w1}\big(\max I^{(1)}_{w}-\min I^{(1)}_{w}
\big)\ge a \big(\max I^{(1)}_{w}-\min I^{(1)}_{w}\big).$$

Since $\max  I^{(2)}> \max I^{(1)}_w$ and $a<1$,
$$\max I^{(2)}-\min  I^{(1)}_{w1}=\max  I^{(2)}-\max I^{(1)}_w+\max
I^{(1)}_w - \min  I^{(1)}_{w1}\ge$$
$$\max  I^{(2)}-\max I^{(1)}_w+ a \big(\max I^{(1)}_{w}
-\min I^{(1)}_{w}\big)\ge$$
$$a\big(\max  I^{(2)}-\max I^{(1)}_w\big)+
a \big(\max I^{(1)}_{w}-\min I^{(1)}_{w}\big)= a\big(\max  I^{(2)}- \min
I^{(1)}_{w}\big),$$
hence
\begin{equation}
\label{eq-eq25}
r(\wii)_{w1}\ge a.
\end{equation}

We next evaluate
$r(\wii)_{w0}$ and $r(\wii)_{w1}$ when (\ref{eq-eq16''}) holds. We have $w=1^{\bar n}$ and
$\max I^{(1)}_{1^{\bar n}}=\max I^{(1)}$. Then,
$$\max I^{(2)}-\min I^{(1)}_w= \max I^{(2)}-\min I^{(2)}+\min I^{(2)}-\max I^{(1)}
+\max I^{(1)}_{1^{\bar n}}- \min I^{(1)}_{1^{\bar n}}$$
$$\le  \max I^{(1)}_{1^{\bar n}}- \min I^{(1)}_{1^{\bar n}}
+{1-a\over a}\big( \max I^{(1)}_{1^{\bar n}}- \min I^{(1)}_{1^{\bar n}}\big)+ \big(
\max I^{(1)}_{1^{\bar n}}- \min I^{(1)}_{1^{\bar n}}\big)$$
$$={a+1\over a} \big( \max I^{(1)}_{1^{\bar n}}- \min I^{(1)}_{1^{\bar
n}}\big)$$

where the inequalities follow from (\ref{eq-eq21}) and (\ref{eq-eq22}). Hence
\begin{equation}
\label{eq-eq26}
r(\wii)_{w0}={\max I_{w0}-\min I_{w0}\over  \max I_{w}-\min I_{w}}=
{\max I^{(1)}_{1^{\bar n}}-\min I^{(1)}_{1^{\bar n}}\over \max I^{(2)}-\min
I^{(1)}_w}\ge {a\over a+1}.
\end{equation}
Similarly,
$$\max I^{(2)}-\min I^{(1)}_w=
\max I^{(2)}-\min I^{(2)}+\min I^{(2)}-\max I^{(1)}
+\max I^{(1)}_{1^{\bar n}}- \min I^{(1)}_{1^{\bar n}}$$
$$\le \max I^{(2)}-\min I^{(2)}+
{1-a\over a} \big(\max I^{(2)}-\min I^{(2)}\big)+ {1\over a}\big(\max I^{(2)}-\min
I^{(2)}\big)$$
$$= {2\over a} \big(\max I^{(2)}-\min I^{(2)}\big).$$
Hence
\begin{equation}
\label{eq-eq27}
r(\wii)_{w1}={\max I_{w1}-\min I_{w1}\over  \max I_{w}-\min I_{w}}=
{\max I^{(2)} -\min I^{(2)} \over \max I^{(2)}-\min I^{(1)}_w}\ge {a\over 2}.
\end{equation}

Finally, we easily see that if (\ref{eq-eq17}), (\ref{eq-eq18}) or (\ref{eq-eq19}) holds, then we have respectively
$r(\wii)_{wj}=r(\wii^{(1)})_{1^{\bar n}w'j}$, $r(\wii)_{wj}=r(\wii^{(2)})_{w'j}$,
  $r(\wii)_{wj}=r(\wii^{(1)})_{ wj}$, and by hypothesis such numbers are
all larger than, or equal to $a$.

To sum up, in view of (\ref{eq-eq24}), (\ref{eq-eq25}), (\ref{eq-eq26}) and (\ref{eq-eq27}), we have $r(\wii)_{wj}\ge a'$ where $a'$ is given by eq. (\ref{eq-aprime})
and the Lemma is completely proven.
\enpr

\section{Finite sums of Cantor Sets}
\label{sec-sum}

In the next Lemma we prove that any finite sum of Cantor sets of ulbd contains a Cantor set of ulbd with the same convex hull as the full sum and maximum gap size not larger than those of the individual Cantor sets.
\begin{lemma}
\label{lem-lem8}
Let $C^{(1)},\ldots,C^{(m)}$ be Cantor sets constructed on $I^{(1)},\ldots,I^{(m)}$
with ratios of dissection larger than, or equal to, $a$. Then there exists a Cantor set
$C\subseteq C^{(1)}+\ldots+C^{(m)}$ constructed with all ratios of
dissection larger than $a_m$, where $a_m>0$ depends only on $a$ and $m$,
such that
\begin{equation}
\gamma(C)\le\max\limits_{s=1,\ldots,m}
\gamma(C^{(s)}) ,
\label{eq-gap}
\end{equation}
\begin{equation}
\min C=\sum\limits_{s=1}^m \min C^{(s)}, \; \max
C=\sum\limits_{s=1}^m \max C^{(s)}.
\label{eq-limits}
\end{equation}
\end{lemma}

\smad
The proof of this Lemma is rather long and technical. For better clarity, it is organized in successive steps.

{\em Step 1.}
The lemma holds trivially for $m=1$. We first show by induction that if it holds for $m=2$, then it holds for any $m$. It will then be sufficient to prove the Lemma for $m=2$. In fact, suppose that the
Lemma holds for $m=2$, and that it also holds for a generic value $m\ge 1$. This implies that it holds for $m+1$, as the following argument shows.
Let $C^{(1)},\ldots,C^{(m)}, C^{(m+1)}$ be Cantor sets
constructed with  all ratios of dissection at least $a$. Then by
hypothesis there exists
a Cantor set $C'\subseteq C^{(1)}+\ldots+C^{(m)}$ which can be constructed with all
ratios of dissection at least $a_m$, such that
$\gamma(C')\le\max\limits_{s=1,\ldots,m} \gamma(C^{(s)})$ and $\min
C'=\sum\limits_{s=1}^m \min C^{(s)}$, $\max C'=\sum\limits_{s=1}^m \max C^{(s)}$.
Put now $a'=\min \{a,a_m\}$ and let $a_{m+1}= a'_2$. Then, by the Lemma for $m=2$ applied to the pair $C'$, $C^{(m+1)}$ there exists a  Cantor set
$$C\subseteq C'+C^{(m+1)}\subseteq   C^{(1)}+\ldots+C^{(m)} +C^{(m+1)}$$
with all dissection ratios at least $a_{m+1}$ such that
$$\gamma(C)\le \max\{\gamma(C'), \gamma(C^{(m+1)})\}\le
\max\limits_{s=1,\ldots,m,m+1} \gamma(C^{(s)}),$$

$$\min C=\min C'+\min C^{(m+1)}
=\sum\limits_{s=1}^{m+1} \min C^{(s)},$$

$$\max C=\max C'+\max C^{(m+1)}
=\sum\limits_{s=1}^{m+1} \max C^{(s)},$$

and the Lemma for $m+1$ holds. In the next steps we will prove the Lemma for $m=2$.

{\em Step 2.}
Let $C^{(1)}, C^{(2)}$ be Cantor sets constructed on $\wii^{(1)}, \wii^{(2)}$
with  all ratios of dissection at least $a$. Based on these latter, we
will define a construction $\wii$ of a new Cantor set $C$. Prior to do that, we need to study auxiliary sets $A_n$.
For any $n \in \bna$ (including obviously $n=0$), put
$$A_n=A^{(1)}_{n}\cup A^{(2)}_{n},$$
where the terms in the union are defined as follows. Suppose that the following condition holds:
\begin{equation}
\gamma(\wii^{(1)})_{0^n}\le \gamma(\wii^{(2)})_{0^n},
\label{eq-eq28}
\end{equation}
$\gamma$ being the gap size defined in eq. (\ref{eq-gapsize}). In this case, let
\begin{equation}
A^{(1)}_{n}:=C^{(2)}_{0^n1}+\max C^{(1)}_{0^{n+1}}, \quad
A^{(2)}_{n}:=   C^{(1)}_{0^n1}    +\max C^{(2)}_{0^n 1}.
\label{eq-eq28b}
\end{equation}
In the opposite case the indices $(1)$ and $(2)$ at r.h.s. of eq. (\ref{eq-eq28b}) are exchanged:
$$A^{(1)}_{n}:= C^{(1)}_{0^n1}+\max C^{(2)}_{0^{n+1}}, \quad
A^{(2)}_{n}:=    C^{(2)}_{0^n1}    +\max C^{(1)}_{0^n 1}.$$
Therefore, let us consider the case in eqs. (\ref{eq-eq28},\ref{eq-eq28b}), the other giving results that can be obtained by exchanging superscripts.
Clearly, $A^{(1)}_{n}$  and $A^{(2)}_{n}$ are two Cantor sets with gaps not larger than
$\max \big\{\gamma(C^{(1)}),  \gamma(C^{(2)}) \big\}$, constructed with ratios
of dissections at least $a$.
We have
$$\min A^{(2)}_n-\max A^{(1)}_n=\min C^{(1)}_{0^n1}    +\max C^{(2)}_{0^n 1}
-\big(\max C^{(2)}_{0^n1}    +\max C^{(1)}_{0^{n+1}}\big)=$$
$$=\min I^{(1)}_{0^n1}+\max I^{(2)}_{0^n 1}-
\max I^{(2)}_{0^n1}-\max I^{(1)}_{0^{n+1}}
 =\min I^{(1)}_{0^n1}-
\max I^{(1)}_{0^{n+1}}=\gamma(\wii^{(1)})_{0^n}.
$$
Since the last quantity is positive, this also proves that $\max A^{(1)}_n < \min A^{(2)}_n$.
Following the same kind of computation,
we also have that
$$\max A^{(2)}_n-\max A^{(1)}_n= \max C^{(1)}_{0^n1}-\max C^{(1)}_{0^n0}= \max I^{(1)}_{0^n1}-\max I^{(1)}_{0^{n}0} \le  $$
$$ \max I^{(1)}_{0^n 1}
- \min I^{(1)}_{0^{n}0} =\max I^{(1)}_{0^n}-\min I^{(1)}_{0^n}=
{1\over r(\wii^{(1)})_{0^n 1}   }
\big(\max I^{(1)}_{0^n 1}-\min I^{(1)}_{0^n 1}\big)=$$
$$={1\over r(\wii^{(1)})_{0^n 1}   }\big(\max A^{(2)}_n-\min A^{(2)}_n\big)\le
{1\over a}\big(\max A^{(2)}_n-\min A^{(2)}_n\big).$$

Therefore,
\begin{equation}
{\max A^{(2)}_n-\min A^{(2)}_n\over \max A^{(2)}_n-\max A^{(1)}_n}\ge a . \label{eq-eq29}
\end{equation}

Moreover
$$\min A^{(2)}_n-\min A^{(1)}_n=\min A^{(2)}_{n}-\max   A^{(1)}_n
+\max A^{(1)}_n-\min A^{(1)}_n=$$
$$\gamma(\wii^{(1)})_{0^n}+
\max C^{(2)}_{0^n1}-\min C^{(2)}_{0^n1}\le$$
$$ \gamma(\wii^{(2)})_{0^n}+
\max I^{(2)}_{0^n1}-\min I^{(2)}_{0^n1}= \min I^{(2)}_{0^n1}-\max
I^{(2)}_{0^{n+1}}+\max I^{(2)}_{0^n1}-\min I^{(2)}_{0^n1} \le$$
$$\max I^{(2)}_{0^n1}-\min I^{(2)}_{0^{n+1}}=
\max I^{(2)}_{0^n}-\min I^{(2)}_{0^n}=
{1\over r(\wii^{(2)})_{0^n 1}   } \big(\max I^{(2)}_{0^n1}-\min
I^{(2)}_{0^n1}\big)$$
$$={1\over r(\wii^{(2)})_{0^n 1}   }\big(\max A^{(1)}_n-\min A^{(1)}_n\big)\le
{1\over a}\big(\max A^{(1)}_n-\min A^{(1)}_n\big).
$$
As a consequence,
\begin{equation}
{\max A^{(1)}_n-\min A^{(1)}_n\over \min A^{(2)}_n-\min A^{(1)}_n}\ge a . \label{eq-eq30}
\end{equation}

Thus, by (\ref{eq-eq29}) and (\ref{eq-eq30}) $A^{(1)}_n$ and $A^{(2)}_n$ satisfy (\ref{eq-eq10}). We can so use Lemma \ref{lem-lem7}, that implies that
$A_n$ is a Cantor set and can be constructed on a map  with all
ratios of dissection larger than, or equal to, $a'$. Let us denote this map with $\wiib(n)$ and its image intervals by $\overline{I}(n)_w$: recall that a different map is defined for any value of $n$, including zero. Moreover, by construction

$$\gamma(A_n)\le \max \big\{ \gamma(A^{(1)}_n), \gamma(A^{(2)}_n),
\min A^{(2)}_n-\max A^{(1)}_n\big\} $$
$$\le  \max \big\{ \gamma(C_1),\gamma(C_2),
\gamma(\wii^{(1)}_{0^n})\big\}\le \max \big\{ \gamma(C_1),\gamma(C_2) \big\}.$$

{\em Step 3.} We can now introduce the new map $\wii$ that constructs the Cantor set $C$ in the thesis of this Lemma. Recall the notation that associates an interval to any finite word, eq (\ref{eq-defwii}): $\wii (w) = I_w$. Let us define all such intervals, parting the set of finite binary words $W$ according to the number of leading zeros. In fact, let
\begin{equation}
\left\{
\begin{array}{l}
I_{0^n}=I^{(1)}_{0^n}+I^{(2)}_{0^n},\\
I_{0^n1w'}=\overline{I}(n)_{w'}.
\end{array}
\right.
\label{eq-defwii+}
\end{equation}
In the above, $w'$ is any word, $n \in \bna$ can take the value 0, and the intervals $\overline{I}(n)_{w'}$, $I^{(1)}_{0^n}$ and $I^{(2)}_{0^n}$ have been defined in the previous step. As before, $0^0$ is to be intended as the empty set.

It is instructive to write down explicitly the first few formulae: let $n=0$, to obtain
$I_\emptyset = I^{(1)}_\emptyset + I^{(2)}_\emptyset$. This is the convex hull of $C$ and it is clearly made by an interval composed of the arithmetic sums of any pair of numbers, one in $I^{(1)}$ and one in $I^{(2)}$. Therefore, it is also the convex hull of $C^{(1)} + C^{(2)}$. Consider next $I_1$. It can be obtained from the second formula in (\ref{eq-defwii+}): $I_1 = \overline{I}(0)_\emptyset$, that is, the convex hull of $A_0$. All intervals corresponding to words starting with 1 are then constructed by the map $\wiib(0)$: in fact, eq. (\ref{eq-defwii+}) yields
$I_{1w} = \overline{I}(0)_w$; as remarked above these intervals construct the Cantor set $A(0)$. Observe that the maximum of this Cantor set is equal to the maximum of $I^{(1)}+I^{(2)}$ and therefore to the maximum of $C^{(1)}+C^{(2)}$. Let us also describe the case $n=1$. This permits to write the interval $I_0$ as $I^{(1)}_0+I^{(2)}_0$. Constructing $A_1$ via the map $\wiib(1)$ then enables us to define all intervals $I_{01w} = \overline{I}(1)_w$, {\em et cetera.}

{\em Step 4.} We now prove formally that eq. (\ref{eq-defwii+}) is a consistent construction of a
quasi Cantor set. Clearly, $I_w$ is an interval for every $w\in W$, so that we just need to prove that for every $w\in W$ we have
\begin{equation}
\min I_w=\min I_{w0}
\label{eq-eq31} 
\end{equation}
\begin{equation}
\max I_{w0}<\min I_{w1}
\label{eq-eq32} 
\end{equation}
\begin{equation}
\max I_{w1}=\max I_w
\label{eq-eq33} 
\end{equation}

If $w=0^n$, we have

  $$\min I_w=\min \big( I^{(1)}_{0^n}+I^{(2)}_{0^n}\big)=
  \min I^{(1)}_{0^n}+\min I^{(2)}_{0^n}$$
  $$=\min I^{(1)}_{0^{n+1}}+\min I^{(2)}_{0^{n+1}}=
  \min \big(I^{(1)}_{0^{n+1}}+I^{(2)}_{0^{n+1}}\big)=\min I_{0^{n+1}}=\min
I_{0^n 0},$$

and (\ref{eq-eq31}) holds. Next, if (\ref{eq-eq28}) holds, then

$$\min I_{w1}-\max I_{w0}=\min\big (C^{(2)}_{0^n1}+\max
C^{(1)}_{0^{n+1}})- \max  \big(I^{(1)}_{0^{n+1}}+I^{(2)}_{0^{n+1}}\big) = $$
$$\min I^{(2)}_{0^n1}+\max I^{(1)}_{0^{n+1}}-\max I^{(1)}_{0^{n+1}}-\max
I^{(2)}_{0^{n+1}}
=\min I^{(2)}_{0^n1}-\max I^{(2)}_{0^{n+1}}=\gamma(\wii^{(2)})_{0^n}$$
and the last quantity is larger than zero. On the contrary, if
(\ref{eq-eq28}) does not hold, we have $\min I_{w1}-\max
I_{w0}=\gamma(\wii^{(1)})_{0^n}>0$, so that in both cases (\ref{eq-eq32}) holds. To sum up, for any word $0^n$ there is $j \in \{1,2\}$ so that
\begin{equation}
 \min I_{0^n1}-\max I_{0^n0}=\gamma(\wii^{(j)})_{0^n}.
\label{eq-eq34} 
\end{equation}
We now prove (\ref{eq-eq33}). We have
$$\max I_w=\max \big( I^{(1)}_{0^n}+I^{(2)}_{0^n}\big)= \max
I^{(1)}_{0^n}+\max I^{(2)}_{0^n}$$
$$=\max I^{(1)}_{0^n 1}+\max I^{(2)}_{0^n 1}=
\max C^{(1)}_{0^n 1}+\max C^{(2)}_{0^n 1}$$
$$=\max A^{(2)}_n=\max A_n=\max \bar{I}(n)=\max I_{0^n 1}=\max I_{w1}$$
and (\ref{eq-eq33}) is proven.

Suppose now $w=0^n1w'$. In this case (\ref{eq-eq31}), (\ref{eq-eq32}) and (\ref{eq-eq33}) follow immediately from the corresponding properties of $\wiib(n)$. In conclusion, the above proves that $\wii$ in fact constructs a quasi Cantor set, which we denote by $C$.

{\em Step 5.} Equation (\ref{eq-eq34}) and a
straightforward argument when $w=0^n1w'$, imply that eq. (\ref{eq-gap}) holds:
\begin{equation}
 \gamma(C)\le\max\{\gamma(C_1),\gamma(C_2)\}.
\label{eq-eq35} 
\end{equation}

Let us now prove eq. (\ref{eq-limits}). We have that
\[
\min C=\min I=\min I_{0^0}=\min  \big(I^{(1)}_{0^0}+
I^{(2)}_{0^0}\big)
=\min I^{(1)}+\min I^{(2)}=\min C^{(1)}+\min C^{(2)}
\]
and similarly
$$\max C=\max  C^{(1)}+\max C^{(2)}.$$
{\em Step 6.}
We now prove that $C\subseteq C^{(1)}+ C^{(2)}$. Take $x\in C$. Then, there exists
an infinite string $\tilde w=i_1 i_2 i_3 \cdots$ such that $x\in I_{i_1 \cdots i_n}$ for
all $n$. We distinguish two cases. If $i_s=0$ for all $s$, then for all $n$,

$$x\in I_{0^n}=I^{(1)}_{0^n}+I^{(2)}_{0^n}=
\big[\min I^{(1)}_{0^n},\max I^{(1)}_{0^n}\big]+ \big[\min I^{(2)}_{0^n},\max
I^{(2)}_{0^n}\big]=$$
$$ [\min I^{(1)}, \min I^{(1)}+\max I^{(1)}_{0^n}-\min I^{(1)}_{0^n}]
+ [\min I^{(2)},\min I^{(2)}+\max I^{(2)}_{0^n}-\min I^{(2)}_{0^n}]$$
$$=\big[\min I^{(1)}+\min I^{(2)}, \min I^{(1)}+\min I^{(2)}+d_n\big]$$

where

$$d_n:=\max I^{(1)}_{0^n}-\min I^{(1)}_{0^n}+\max I^{(2)}_{0^n}-\min
I^{(2)}_{0^n} \llim\limits_{n\to +\infty}0,$$
since $\wii^{(1)}$ and $\wii^{(2)}$ construct two Cantor sets.
Hence, $x=\min I^{(1)}+\min I^{(2)}=\min  C^{(1)}+\min C^{(2)}\in
  C^{(1)}+ C^{(2)}$. Suppose instead, there exists $\bar s$ such that
$i_{\bar s}=1$, and we can and do assume $i_s=0$ for every $s<\bar s$ (in other words, ${\bar s}$ is the first occurrence of 1 in the symbolic sequence $\tilde w$). Putting
$n=\bar s-1$ and $\tilde w'= i_{\bar s+1}i_{\bar s+2} \cdots$, we have $\tilde w_{\bar
s+m}=0^n1\tilde w'_{m}$
  for all $m>0$. Thus for every $m>0$,

  $$x\in I_{0^n1\tilde w'_m}=\overline{I}(n)_{\tilde w'_m}\Rightarrow x\in
A_n=A^{(1)}_n\cup
  A^{(2)}_n\subseteq C^{(1)}+C^{(2)}. $$

Thus, $C\subseteq C^{(1)}+C^{(2)}$ is proven.

{\em Step 7.}
It finally remains to prove that
$r(\wii)_w\ge a'$ for all $w\in W$ (as usual this also proves that $C$ is not only
quasi Cantor but also Cantor).
Let us start by considering the word $w=0^n$. Recalling that $r(\wii^{(j)})_v\ge a$,
for $j=1,2$ and $v\in W\setminus\{\emptyset\}$, we have the following estimates.
The diameter of the interval $I_w$ is:
\begin{eqnarray}
 d(I_w)=d(I_{0^n})= \max \big(I^{(1)}_{0^n}+I^{(2)}_{0^n}\big)-\min
\big(I^{(1)}_{0^n}+I^{(2)}_{0^n}\big)= \\
\max I^{(1)}_{0^n}-\min I^{(1)}_{0^n}+
\max I^{(2)}_{0^n}-\min I^{(2)}_{0^n}=d(I^{(1)}_w)+d(I^{(2)}_w).
\end{eqnarray}
The diameter of the interval $I_{w0}$ is:
\begin{eqnarray}
d(I_{w0})=d(I_{0^{n+1}})=d(I^{(1)}_{w0})+d(I^{(2)}_{w0})
=r(\wii^{(1)})_{w0}d(I^{(1)}_w)+
r(\wii^{(2)})_{w0}d(I^{(2)}_w) \\
\ge a d(I^{(1)}_w)+ a d(I^{(2)}_w)=a d(I_w)
\end{eqnarray}
The diameter of the interval $I_{w1}$ is:
\begin{eqnarray}
d(I_{w1})=d(I_{0^n1})=d\big(\overline{I}(n)\big)=\max \overline{I}(n)-\min \overline{I}(n)= \\
=\max A_n-\min A_n=\max A^{(2)}_n-\min A^{(1)}_n \\
=\max C^{(1)}_{0^n1}+\max C^{(2)}_{0^n1}
-\min C^{(2)}_{0^n1}-\max C^{(1)}_{0^{n+1}}.
\end{eqnarray}
If  (\ref{eq-eq28}) holds, we can continue as follows:
\begin{eqnarray}
d(I_{w1}) =\max I^{(2)}_{0^n1}-\min I^{(2)}_{0^n1}+\max I^{(1)}_{0^n1}
-\max I^{(1)}_{0^{n+1}} \\
\ge \max I^{(2)}_{0^n1}-\min I^{(2)}_{0^n1}+\max I^{(1)}_{0^n1}-\min
I^{(1)}_{0^n1} \\
=d(I^{(2)}_{0^n1})+ d(I^{(1)}_{0^n1})=d(I^{(2)}_{w1})+ d(I^{(1)}_{w1}) \\
=r(\wii^{(2)})_{w1} d(I^{(2)}_w)+ r(\wii^{(1)})_{w1} d(I^{(1)}_w) \\
\ge
 a d(I^{(2)}_w)+ a d(I^{(1)}_w)= a d(I_w).
 \end{eqnarray}

Hence,

$$r(\wii)_{wj}={d(I_{wj})\over d(I_w)}\ge a\ge a', \quad
j=0,1.$$

Note that for $j=1$ we have used eq. (\ref{eq-eq28}); when it does not hold we exchange the indices $1$ and $2$.

Let now consider the words $w=0^n1w'$, with $w'$ any finite word. Recalling that
$r\big(\wiib(n)\big)_v\ge a'$, for $n\in\bna$ and $v\in W\setminus\{\emptyset\}$, we
have
\[r(\wii)_{wj}={d(\overline{I}_{wj})\over d(\overline{I}_w)}=
{d\big(\overline{I}(n)_{w'j}\big)\over d\big(\overline{I}(n)_{w'}\big)}=r\big(\wiib(n)\big)_{w'j}\ge a'
\]
and the Lemma is completely proven.  \enpr

\section{More on sums of Cantor sets}
\label{sec-more}

In this section, we recall a result from the literature on the arithmetical sum of Cantor sets, and we use it to derive a couple of further Lemmas.

\begin{theorem}[Theorem 3.2 in \cite{cabrelli}]
\label{lem-lem5}
Let $C^{(1)},\ldots,C^{(m)}$ be Cantor sets constructed on $I^{(1)},\ldots,I^{(m)}$
with constructions $\wii^{(1)},\ldots,\wii^{(m)}$  of uniformly lower bounded dissections, larger than $a>0$.
In addition, suppose that $a\le {1\over 3}$ and that $m$ is
such that $\disp{(m-1){a^2\over (1-a)^3}+ {a\over 1-a}\ge 1}$. Finally suppose
that no translate of any of these Cantor sets is contained in a gap of another. Then, the sum of these Cantor sets is a closed interval:
\[
C^{(1)}+ \ldots +C^{(m)}=\big[\sum\limits_{i=1}^m \min C^{(i)},
    \sum\limits_{i=1}^m \max C^{(i)}\big].
\]
\end{theorem}

The result of this theorem can be easily extended to any sets containing the Cantor sets $C^{(i)}$ and enclosed in the intervals $I^{(i)}$.

\begin{lemma}
\label{lem-lem6}
Let $C^{(i)}$ for $i=1,\ldots,m$, be Cantor sets and suppose that $C^{(1)}+ \ldots +C^{(m)}$ is an interval.
Then, for any sets $D^{(i)}$ such that
\begin{equation}
\label{eq-eq7}
C^{(i)}\subseteq D^{(i)}\subseteq I^{(i)},\ \ \    i=1,\ldots,m,
\end{equation}
and
\begin{equation}
\label{eq-eq8}
\min C^{(i)}=\min D^{(i)}, \max C^{(i)}=\max D^{(i)},\ \ \   i=1,\ldots,m,
\end{equation}
we have that
\begin{equation}
\label{eq-eq9}
D^{(1)}+\ldots.+D^{(m)}=\big[\sum\limits_{i=1}^m \min C^{(i)}, \sum\limits_{i=1}^m \max
C^{(i)}\big].
\end{equation}
\end{lemma}

\smad Observe that, by eq. (\ref{eq-eq7}),

$$C^{(1)}+\ldots.+C^{(m)}\subseteq D^{(1)}+\ldots+D^{(m)}\subseteq I^{(1)}+\ldots+I^{(m)},$$

Moreover, we trivially have
$$I_1+\ldots+I_m=\big[\sum\limits_{i=1}^m \min I^{(i)}, \sum\limits_{i=1}^m \max
I^{(i)}\big]= \big[\sum\limits_{i=1}^m \min C^{(i)}, \sum\limits_{i=1}^m \max C^{(i)}\big]$$

and, since $C^{(1)}+\ldots+C^{(m)}$ is an interval, it is clearly
$$C^{(1)}+\ldots+C^{(m)}=\big[\sum\limits_{i=1}^m \min C^{(i)}, \sum\limits_{i=1}^m \max
C^{(i)}\big].  \enpr$$

We now use Theorem \ref{lem-lem5} to prove the following:

\begin{lemma}
\label{lem-lem9}
Suppose that $C^{(l)}$, $l=1,2,\ldots$ are Cantor sets of ulbd larger than $a>0$, such that
$d\big(C^{(l)}\big)\in [A_1,A_2]$ for any $l$, with  $0<A_1< A_2<\infty$. Then, there exists $n \in \bna$, that depends only on $A_1,A_2$ and $a$, such that
$$C^{(1)}+ C^{(2)}+\cdots+C^{(n)}=
\Big[\sum\limits_{l=1}^n \min C^{(l)}, \sum\limits_{l=1}^n \max C^{(l)}\Big].$$
\end{lemma}

\smad Take the smallest $m \in \bna$ such that  $mA_1>A_2$. Let $h \in \bna$ and consider the finite sums $S^{(h)}= C^{(hm+1)}+\cdots+ C^{(hm+m)}$. Because of
Lemma \ref{lem-lem8}, for every $h\in\bna$
there exist a Cantor set $D^{(h)}$, constructed with ratios of dissection larger than, or equal to $a_m$ (where $a_m$ do not depend on $h$ and where, of course, we can take $a_m\in ]0,{1\over 3}]$) such that
\begin{equation}
\label{eq-eq36}
  D^{(h)}\subseteq C^{(hm+1)}+\cdots+ C^{(hm+m)}
\end{equation}

\begin{equation}
\label{eq-eq37}
 \min D^{(h)}=\sum\limits_{l=hm+1}^{hm+m} \min C^{(l)},
\end{equation}

\begin{equation}
\label{eq-eq38}
 \max D^{(h)}=\sum\limits_{l=hm+1}^{hm+m} \max C^{(l)},
\end{equation}

\begin{equation}
\label{eq-eq39}
 \gamma(D^{(h)})\le\max\limits_{l=hm+1,\ldots,hm+m} \gamma(C^{(l)})\le A_2.
\end{equation}

By (\ref{eq-eq37}) and (\ref{eq-eq38}) we have
$$d(D^{(h)})=\max D^{(h)}-\min D^{(h)}
=\sum\limits_{l=hm+1}^{hm+m} d(C^{(l)})\ge mA_1>A_2,
$$

so that, in view of (\ref{eq-eq39}) $d(D^{(h)}) > \gamma(D^{(h')})$ for any $h,h' \in \bna$, {\em i.e.} no translate of $D^{(h)}$ is contained in a gap of
$D^{(h')}$. Now, by Theorem \ref{lem-lem5}, if $H \in \bna$ is large enough, then

$$D^{(0)}+\cdots+D^{(H)}=
\Big[ \sum\limits_{h=0}^H\min D^{(h)}, \sum\limits_{h=0}^H\max D^{(h)}\Big],$$

so that on one hand

$$C^{(1)}+\cdots + C^{(Hm+m)}\subseteq
\Big[\sum\limits_{i=1}^{Hm+m}\min C^{(i)}, \sum\limits_{i=1}^{Hm+m}\max
C^{(i)}\Big].$$

On the other hand, by (\ref{eq-eq36}) and (\ref{eq-eq37})
$$C^{(1)}+\cdots + C^{(Hm+m)}\supseteq D^{(0)}+\cdots+D^{(H)}=
\Big[ \sum\limits_{h=0}^H\min D^{(h)}, \sum\limits_{h=0}^H\max D^{(h)}\Big]$$
$$=\Big[\sum\limits_{i=1}^{Hm+m}\min C^{(i)},
\sum\limits_{i=1}^{Hm+m}\max C^{(i)}\Big]
$$
and the Lemma is proven for $n=Hm+m$. \enpr

\section{Proof of the main theorem and related results}
\label{sec-proof}

In this section we prove the main theorem, via an additional proposition, where we describe a situation where the geometric sum of Cantor sets is an interval. Since in some sense this is a generalization of the results in \cite{cabrelli}, we think it is interesting in itself.

\begin{proposition}
\label{lem-lem10}
Let $K$ be a ulbd Cantor set of dissection ratios larger than $a>0$, or the attractor of a finite, hyperbolic non--singular I.F.S. with maps (\ref{eq-contra}). Then, for any $0<\alpha<1$ there exists $n \in \bna$  such that
$K_n={1\over\alpha}K+\ldots+{1\over\alpha^n}K$ is the disjoint union of finitely many closed intervals.
\end{proposition}

\smad

Let us consider the first case: let $K$ be a ulbd Cantor set. Then, it can be constructed by a map $\wii$ on an initial interval $I_\emptyset$. Let $A=d(I_\emptyset)$ and consider the following condition: for a given $l \in \bna$  require that
\begin{equation}
\label{eq-cond1}
 d(I_w) \geq A \alpha^l,  \; d(I_{wj}) < A \alpha^l
\end{equation}
for at least one value of $j \in \{0,1\}$.  Clearly, this condition may or may not be verified by a finite word $w$. Let $\ccb_l$ the set of words that pass the test, for a given value $l \in \bna$:
\begin{equation}
\label{eq-cond2}
  \ccb_l = \{ w \in W \mbox{  s.t. (\ref{eq-cond1}) holds} \}
\end{equation}
Recall that $\wii_w$ is the induced construction on $I_w$ defined in eq. (\ref{eq-defwii2}) yielding the Cantor set $C(\wii_w)$. Then, it is not difficult to see that for any choice of $l \in \bna$, $K$ may be written as the union of a finite number of Cantor sets, as follows:
\begin{equation}
\label{eq-cond3}
 K = \bigcup_{w \in \ccb_l}  C(\wii_w).
\end{equation}
In fact, because of eq. (\ref{eq-eq4x}), for any $x \in K$, for any $n \in \bna$ there exists a word $w_n(x)$ of length $n$ such that $x \in I_{w_n(x)}$, hence $x\in C\big({\cal I}_{w_n(x)}\big)$ (see eq. (\ref{eq-ciccia}). Since $d(I_{w_n(x)})$ tends to zero when $n$ tends to infinity, there exists at least one
value $n$ such that (\ref{eq-cond1}) holds. This proves inclusion. Furthermore, since zero is the unique accumulation points of diameters $d(I_{w})$ the cardinality of $\ccb_l$ is finite.

Let $w \in \ccb_l$. Condition (\ref{eq-cond1}) implies that $d(I_w) \geq A \alpha^l$. Moreover, there exists $j \in \{0,1\}$ such that
 $  A \alpha^l > d(I_{wj}) = r(\wii)_{wj} d(I_w) \geq a d(I_w), $
so that
\begin{equation}
\label{eq-cond4b}
  A \alpha^l \leq  d(C(\wii_w)) = d(I_w) < \frac{A}{a} \alpha^l .
\end{equation}
Next, divide all terms in eq. (\ref{eq-cond3}) by $\alpha^l$, to prove that $K/\alpha^l$ can be written as a finite union of Cantor sets $C(\wii_w)/\alpha^l$, $w \in \ccb_l$, each of which has the following properties: it has uniform lower bounded dissection larger that $a>0$, and its diameter lies in the interval $[A,A/a]$. Observe that this holds for any value $l \in \bna$.

Consider now the set $K_n = {1\over\alpha}K+\ldots+{1\over\alpha^n}K $ in the thesis of the Lemma. Using eq. (\ref{eq-cond3}) it can be written as follows:
\begin{equation}
\label{eq-cond5}
 K_n = \bigcup_{w_1 \in \ccb_1, \ldots, w_n \in \ccb_n} \frac{C(\wii_{w_1})}{\alpha} + \ldots +  \frac{C(\wii_{w_n})}{\alpha^n}.
\end{equation}
We can now apply Lemma \ref{lem-lem9} to prove that there exists an integer $n$ that depends only on $a$ and $A$ such that $\frac{C(\wii_{w_0})}{\alpha} + \ldots +  \frac{C(\wii_{w_n})}{\alpha^n}$ is an interval. Since the cardinality of each $\ccb_l$ is finite, it follows at once that $K_n$ is the union of a finite number of disjoint, closed intervals.

More complicated is the case when $K=K_\Psi$ is the attractor of a non--singular finite hyperbolic I.F.S. Denote again by $I_\emptyset$ the convex hull of $K$, $A=d(K)$. Let $\delta$ and $\sigma$ be as in eq. (\ref{eq-contra}).
Also denote by $V$ the set of finite words composed of the letters $\{1,\ldots,M\}$ and by $\psi_v$ the composite map defined in eq. (\ref{eq-map1x}) with $v=i_1,\cdots,i_n$. Consider the diameters of the sets $\psi_v(K)$. Clearly, $d(\psi_v(K))=d(\psi_v(I_\emptyset))$: to these latter we can apply  Lemma \ref{lem-issimo}, which
proves that letting $c=\exp(-C)$, $C>0$, computed as in eq. (\ref{eq-crat6}) we have
\begin{equation}
   d(\psi_{vi}(K)) \geq  c d(\psi_{v}(K))
\label{eq-cicci}
\end{equation}
for any $v \in V$, $i=1,\ldots,M$, being $vi$ the composed word.
We can now replace condition (\ref{eq-cond1}) by the following: for $l \in \bna \setminus \{0\}$ require that $v \in V$ satisfies
 \begin{equation}
\label{eq-cond1bis}
 c A \alpha^l < d(\psi_v(K)) \leq A \alpha^l,
\end{equation}
and define accordingly
\begin{equation}
\label{eq-cond2bis}
  \ccb_l = \{ v \in V \mbox{  s.t. eq. (\ref{eq-cond1bis}) holds} \}.
\end{equation}
The analogue of eq. (\ref{eq-cond3}) is now
\begin{equation}
\label{eq-cond3bis}
 K = \bigcup_{v \in \ccb_l}  \psi_v(K).
\end{equation}
To prove eq. (\ref{eq-cond3bis}) observe that for any $x\in K$ there exists $i_1$ such that $x\in \psi_{i_1}(K)$, then
there exists $i_2$ such that $x\in  \psi_{i_1,i_2}(K)$ and so proceeding: for every
$s$ there exist $i_1,\ldots,i_s$ such that $x\in \psi_{i_1,\ldots,i_s}(K)$. Choose the first index $s>0$ such that $d(\psi_{i_1,\ldots,i_s}(K))\le \alpha^l d(K)$, which surely exists, since $d(\psi_{i_1,\ldots,i_s}(K))\le \delta^s d(K)$.
At the same time, since $d(\psi_{i_1,\ldots,i_{s-1}}(K)) > \alpha^l d(K)$ eq. (\ref{eq-cicci}) implies that
\[
d(\psi_{i_1,\ldots,i_s}(K)) \geq c d(\psi_{i_1,\ldots,i_{s-1}}(K)) > c \alpha^l d(K)
\]
so that (\ref{eq-cond1bis}) holds for $v=i_1,\ldots,i_s$ and therefore $K$ is a subset of the union at r.h.s. of (\ref{eq-cond3bis}). The other inclusion is obvious from $K = U_\Psi^n (K)$.
The fact that $\ccb_l$ has only finitely many elements follows easily from the inequality $d(\psi_{i_1,\ldots,i_s}(K)) \le \delta^s d(K)$ and from the first inequality in condition (\ref{eq-cond1bis}).

We can now write the analogue of eq. (\ref{eq-cond5}):
\begin{equation}
\label{eq-cond5bis}
 K_n = \bigcup_{v_1 \in \ccb_1, \ldots, v_n \in \ccb_n} \frac{\psi_{v_1}(K)}{\alpha} + \ldots +  \frac{\psi_{v_n}(K)}{\alpha^n}.
\end{equation}
We need to analyze this union. At difference with the first part of the proof, we cannot apply Lemma \ref{lem-lem9} directly, because we have control of the diameter of the sets ${\psi_{v_l}(K)}/{\alpha^l}$ (that is contained in the interval $[cA,A]$) but not of their nature.

Let $\Psi'$ the two-maps I.F.S. related to $\Psi$ as in Lemma \ref{lem-3}. Recall that $\Psi'$ is obtained by selecting two maps out of the full set $\Psi$, in such a way to conserve the convex hull of the attractor. By a suitable relabeling of indices let $\Psi' = \{\psi_0,\psi_1\}$. Let $K' = K_{\Psi'}$ be its attractor. We have $K'\subseteq K$ (Lemma \ref{lem-1} iii).
$K'$ is either a closed interval, or a Cantor set, with the same convex hull of $K$. In the former case $K=K'$ is the same closed interval, and therefore ${1\over\alpha}K+\ldots+{1\over\alpha^n}K$ is an interval, and the thesis of this Lemma follows easily.

The second case is more interesting. $\Psi'$ is a two--maps, non--singular I.F.S., whose attractor $K'$ is a Cantor set. Lemma \ref{lem-nuovo} establishes that $K'$ has a ulbd construction. Consider now the images
$K'_v:=\psi_v(K')$, with $v \in V$. Each of these is a Cantor set. We can easily prove that they too are of ulbd. In fact, a construction $\wii^v$ for each of them can be obtained from $\wii$ in Lemma \ref{lem-nuovo}, in analogy with eq. (\ref{eq-defwii2}) as:
\[
   I^v_{h(w)} = \psi_v(I_w) = (\psi_v \circ \psi_w)(I_\emptyset) = \psi_{vw} (I_\emptyset),
\]
where $w$ is a finite word in the labels of the two maps that compose $\Psi'$ and $h$ is a permutation of the finite word $w$, constructed along the same lines of lemma \ref{lem-nuovo}.
Lemma \ref{lem-issimo}, which we have also used above, proves that the dissection ratios of $K'_v:=\psi_v(K')$ are uniformly lower bounded by the value $c = \exp(-C) >0$, eq. (\ref{eq-crat6}), computed over the full set of maps composing $\Psi$.

Let us now replace $K$ by $K'$ at r.h.s. in eq. (\ref{eq-cond5bis}).
Since $K$ and $K'$ have the same convex hull, conditions (\ref{eq-cond1bis}) and (\ref{eq-cond2bis}) imply that the diameters of ${\psi_{w_l}(K')}/{\alpha^l}$ are all contained in $[c A,A]$. Lemma \ref{lem-lem9} then implies that for sufficiently large $n$ (that depends only on $c$, $\delta$ and $A$) the finite sum
${1\over\alpha} \psi_{v_1}(K')+\ldots+{1\over\alpha^n} \psi_{v_n}(K')$ is an interval, for any choice of $v_1,\ldots,v_n$ in the respective sets. We can now apply Lemma \ref{lem-lem6}: setting $C^{(i)} = {1\over\alpha^i} \psi_{v_i}(K')$, $D^{(i)} = {1\over\alpha^i} \psi_{v_i}(K)$, $I_\emptyset$ the convex hull of $K$ and $I^{(i)} = {1\over\alpha^i} \psi_{v_i}(I_\emptyset)$ we are in the conditions of the Lemma, and we obtain the equality
\[
 {1\over\alpha} \psi_{v_1}(K)+\ldots+{1\over\alpha^n} \psi_{v_n}(K) =
 {1\over\alpha} \psi_{v_1}(K')+\ldots+{1\over\alpha^n} \psi_{v_n}(K'),
 \]
which shows that the l.h.s. is also an interval.  Then the Lemma follows again by the finite cardinality of each $\ccb_l$. \enpr

We can now prove the main result of this paper, Theorem
\ref{teo-1}.

Let  $K_\Phi$ be the attractor of the I.F.S. composed of the maps $\Phi$ in eqs. (\ref{eq-map1b}), which is the unique solution in $\ck$ of the equation $\bfuf(K_\Phi)=K_\Phi$.
Note that of course $\bfuf^n(K_\Phi)=K_\Phi$, where the map $\bfuf^n$ is defined by
\begin{equation}
\label{eq-eq4}
\bfuf^n(A)
={
\bigcup\limits_{\beta_n,\ldots,\beta_1\in K}
\phi_{\beta_n}\circ\cdots\circ \phi_{\beta_1}(A)
}
\quad\forall\, A\in \ck.
\end{equation}

Equation (\ref{eq-eq4}) easily follows by induction on $n$.
We now need an algebraic formula: we therefore interrupt the proof of the theorem for the last Lemma:

\begin{lemma}
\label{lem-4}
The n-fold map composition in eq. (\ref{eq-eq4}) takes the following form:
for every $x\in\bre$
\begin{equation}
\label{eq-algebra}
\phi_{\beta_n}\circ\cdots\circ \phi_{\beta_1}(x)=\alpha^n x+
\alpha^{n} (1-\alpha)\sum\limits_{i=1}^n{\beta_i\over \alpha^i} .
\end{equation}
\end{lemma}

\smad We proceed by induction. For $n=0$ and $n=1$ the result is trivial. Suppose
it holds for $n$. For $n+1$ we have

$$\phi_{\beta_{n+1}}\circ\phi_{\beta_n}\circ\cdots\circ \phi_{\beta_1}(x)=
\phi_{\beta_{n+1}}\Big(\phi_{\beta_n}\circ\cdots\circ \phi_{\beta_1}(x)\Big)=$$
$$\phi_{\beta_{n+1}}\Big(\alpha^n x+
\alpha^{n} (1-\alpha)\sum\limits_{i=1}^n{\beta_i\over \alpha^i}\Big)=
\alpha \Big(\alpha^n x+
\alpha^{n} (1-\alpha)\sum\limits_{i=1}^n{\beta_i\over \alpha^i}\Big)+(1-\alpha)
\beta_{n+1}=$$
$$\alpha^{n+1} x+
\alpha^{n+1} (1-\alpha)\sum\limits_{i=1}^{n+1}{\beta_i\over \alpha^i}\, . \enpr
$$

{\em Continuation of the proof of the theorem}. Because of eq. (\ref{eq-eq4}),

\begin{equation}
\label{eq-eq4b}
K_\Phi
={
\bigcup\limits_{\beta_n,\ldots,\beta_1\in K}
\phi_{\beta_n}\circ\cdots\circ \phi_{\beta_1}(K_\Phi) =
\bigcup\limits_{x \in K_\Phi}
\bigcup\limits_{\beta_n,\ldots,\beta_1\in K} \{
\phi_{\beta_n}\circ\cdots\circ \phi_{\beta_1}(x) \}.
}
\end{equation}
Use now eq. (\ref{eq-algebra}) to get
\begin{equation}
\label{eq-eq4c}
\bigcup\limits_{\beta_n,\ldots,\beta_1\in K}
\{
\phi_{\beta_n}\circ\cdots\circ \phi_{\beta_1}(x) \} =
\alpha^n x+
\alpha^{n} (1-\alpha)\sum\limits_{i=1}^n{K\over \alpha^i}.
\end{equation}
Because of Proposition \ref{lem-lem10} we can choose $n$ in such a way that
$\alpha^{n} (1-\alpha)\sum\limits_{i=1}^n{K\over \alpha^i}$ is a finite union of closed intervals, call it $\cal J$. 
The term $\alpha^n x$ in eq. (\ref{eq-eq4c}) merely shifts these intervals. We can now go back to eq. (\ref{eq-eq4b}): since $K_\Phi$ is bounded and closed, we read that $K_\Phi$ is the union of shifted intervals in $\cal J$, and therefore it is itself a finite union of closed intervals.
\enpr

\end{document}